\documentclass[a4paper,10pt,oneside,fleqn]{article}

\usepackage[T1]{fontenc}
\usepackage[a4paper, left=3cm, right=3cm, top=3cm, bottom=3cm, ]{geometry}
\usepackage[english]{babel}
\usepackage{lmodern}

\usepackage{graphicx}
\usepackage{amsmath,amsfonts,amssymb}
\usepackage{mathabx} 
\usepackage{graphics}
\usepackage{fancyhdr}
\usepackage{setspace}
\usepackage{lastpage}
\usepackage[colorlinks,linkcolor=blue,citecolor=blue,urlcolor=blue]{hyperref}

\usepackage{titlesec}
\titleformat*{\section}{\large\bfseries}
\titleformat*{\subsection}{\normalsize\bfseries}
\titlespacing*{\section}{0pt}{3ex plus 1ex minus .2ex}{0ex plus .2ex}
\titlespacing*{\subsection}{0pt}{1ex plus 0ex minus .2ex}{0ex minus 0.5ex}
\titlespacing*{\subsubsection}{0pt}{1ex plus 0ex minus .2ex}{0ex minus 1ex}

\usepackage{tabularx}
\usepackage{bm}
\usepackage{color}
\usepackage{units}
\usepackage{float}
\usepackage{xcolor}
\usepackage{tikz}
\usetikzlibrary{shapes.geometric,arrows}
\usetikzlibrary{calc}
\usetikzlibrary{positioning}

\usepackage{units}
\usepackage{datetime}
\usepackage{booktabs}
\usepackage{array}
\usepackage{multirow}
\usepackage{cite}
\usepackage{caption}
\usepackage{subcaption}

\usepackage{cleveref}
\crefname{figure}{Fig.}{Figs.}
\crefname{table}{Table}{Tables}
\crefname{equation}{Eq.}{Eqs.}
\Crefname{equation}{Equation}{Equations}
\crefname{chapter}{Chapter}{Chapters}
\crefname{appendix}{Appendix}{Appendices}
\crefname{section}{Section}{Sections}
\crefname{subsection}{Section}{Sections}
\crefname{subsubsection}{Section}{Sections}

\newcommand{\ddt}[1]{\frac{\mathrm{d}#1}{\mathrm{d}t}}
\newcommand{\ind}[1]{_{\mathrm{#1}}}
\newcommand{\etal}[0]{et al.~}
\newcommand{\lie}{\mathcal{L}}

\newcounter{savefootnote}
\newcounter{symfootnote}
\newcommand{\symfootnote}[1]{%
	\setcounter{savefootnote}{\value{footnote}}%
	\setcounter{footnote}{\value{symfootnote}}%
	\ifnum\value{footnote}>8\setcounter{footnote}{0}\fi%
	\let\oldthefootnote=\thefootnote%
	\renewcommand{\thefootnote}{\fnsymbol{footnote}}%
	\footnote{#1}%
	\let\thefootnote=\oldthefootnote%
	\setcounter{symfootnote}{\value{footnote}}%
	\setcounter{footnote}{\value{savefootnote}}%
}

\addto\captionsenglish{}

\newenvironment{acknowledgements}{\noindent\footnotesize\textbf{Acknowledgements}}{}
\newenvironment{keywords}{\noindent\textbf{Keywords:}}{}
\renewenvironment{abstract}{\noindent\textbf{Abstract:}}{}

\begin{document}

\begin{center}
	\begin{Large}
		\textbf{Identification of MIMO Wiener-type Koopman Models for Data-Driven Model Reduction using Deep Learning}
	\end{Large} 
\end{center}

\begin{center}
	\begin{large}
		Jan C.~Schulze$^1$, Danimir T.~Doncevic$^2$ and Alexander Mitsos$^{1,3,2,}\symfootnote{Corresponding author: amitsos@alum.mit.edu}$
	\end{large}
\end{center}   

\begin{flushleft}
	\begin{small}
		\begin{tabular}{l p{13.5cm} }
			$^1$ & Process Systems Engineering (AVT.SVT), RWTH Aachen University, 52074 Aachen, Germany.\\
			$^2$ & Institute of Energy and Climate Research: Energy Systems Engineering (IEK-10), Forschungszentrum J\"ulich GmbH, 52425 J\"ulich, Germany.\\
			$^3$ & JARA-ENERGY, 52056 Aachen, Germany.\\
		\end{tabular}
	\end{small}
\end{flushleft}

\vspace{2ex}
\begin{abstract}
We use Koopman theory to develop a data-driven nonlinear model reduction and identification strategy for multiple-input multiple-output (MIMO) input-affine dynamical systems.
While the present literature has focused on linear and bilinear Koopman models, we derive and use a Wiener-type Koopman formulation.
We discuss that the Wiener structure is particularly suitable for model reduction, and can be naturally derived from Koopman theory.
Moreover, the Wiener block-structure
unifies the mathematical simplicity of linear dynamical blocks and the accuracy of bilinear dynamics.
We present a Koopman deep-learning strategy combining autoencoders and linear dynamics that generates low-order surrogate models of MIMO Wiener type. 
In three case studies, we apply our framework for identification and reduction of a system with input multiplicity, a chemical reactor and a high-purity distillation column.
We compare the prediction performance of the identified Wiener models to linear and bilinear Koopman models.
We observe the highest accuracy and strongest model reduction capabilities of low-order Wiener-type Koopman models, making them promising for control.
\end{abstract}

\vspace{2ex}

{
	\footnotesize
	Published in \textit{Computers and Chemical Engineering}:
	\href{https://doi.org/10.1016/j.compchemeng.2022.107781}{DOI:10.1016/j.compchemeng.2022.107781}\\
	\textcopyright{} 2022. This manuscript version is made available under the CC-BY-NC-ND 4.0 license (\url{https://creativecommons.org/licenses/by-nc-nd/4.0/}).\\		
}

\section{Introduction}
The identification of low-order nonlinear dynamical models is one of the main challenges in control engineering.
Practical control models do not have to be highly accurate and broadly applicable, but should be simple and reliable predictors of the control-relevant dynamical response.
Traditionally, low-order control models are obtained either by combining mechanistic modeling and model reduction \cite{Marquardt.2002, Antoulas.2005b}, or via black-box identification of generic input-output models from measured data \cite{Ljung.1999, Pearson.2003, Schoukens.2017}.
Nonlinear model reduction remains an active field of research, as no universal approach exists and the application of many existing methods can be mathematically involved \cite{Antoulas.2005b, Baur.2014}.
On the other hand, black-box approaches such as NARMAX, Volterra series, and Hammerstein or Wiener models 
are nowadays established methods
\cite{Billings.2013,Schoukens.2019}.
However, successful identification is strongly system-specific as it depends on a suitable choice of model type, data set, algorithm as well as user experience \cite{Pearson.2003}.
Recently, data-driven strategies combing elements of the two traditional approaches have been proposed, e.g., dynamic mode decomposition \cite{Schmid.2010, Williams.2015, Klus.2018}, Koopman-based modeling \cite{Mauroy.2020}, lift-and-learn methods \cite{Qian.2020},  and deep-learning autoencoder approaches \cite{Lee.2020, Masti.2021}.
Overall, generating low-order control models remains challenging in practice due to difficulties in the specification and interpretation of reduced model structures, and involved reduction or identification algorithms.

Here, we perform nonlinear model reduction through data-driven identification of low-order surrogate models.
For this purpose, we apply Koopman theory and develop a broadly applicable reduction strategy that identifies low-order MIMO Wiener control models.
These reduced models are simple and relatively minimalistic due to their linear dynamical block.
Our approach employs deep-learning methods and targets a minimal need of expert knowledge of model reduction.

Koopman theory postulates that nonlinear dynamical systems can be
globally exactly reformulated as linear dynamics
when represented on an infinite-dimensional state space by means of nonlinear coordinate transformation \cite{Mezic.2005, Otto.2021}.
Recent works have aimed at identifying finite-dimensional, possibly low-order approximates of the infinite-dimensional Koopman operator \cite{Brunton.2016b}.
Further, Koopman theory has been adopted for nonlinear system analysis, e.g., global stability analysis \cite{Mauroy.2016} and the extension of the local Hartman-Grobman theorem to the entire basin of attraction of an equilibrium point or periodic orbit \cite{Lan.2013}. 
While finding an exact low-dimensional Koopman model is often difficult or impossible,
using data-driven methods, e.g., deep learning, to identify approximate Koopman models for monitoring or control is feasible in many cases.
Examples for the successful applications of Koopman theory in various fields include fluid dynamics \cite{Rowley.2009},
molecular thermodynamics \cite{Wehmeyer.2018} 
and disease modeling \cite{Proctor.2015}.
Due to its recent extension to systems with external inputs \cite{Surana.2016, Goswami.2017, Proctor.2018, Villanueva.2021}, 
Koopman theory is also promising for control.

We present a data-driven non-intrusive reduction strategy for control motivated by Koopman theory. 
Therein, we simultaneously identify a nonlinear transformation and low-order latent dynamics approximating the nonlinear system evolution.
While the existing literature on Koopman modeling for control has focused on 
linear \cite{Korda.2018, Arbabi.2018, Narasingam.2019, Abraham.2019, Folkestad.2020, Han.2020, Son.2021b, Ping.2021, Lin.2021} and bilinear \cite{Peitz.2020, Bruder.2021, Folkestad.2021, Kaiser.2021} forms, 
we consider a MIMO Wiener structure.
As we show, when combined with Koopman theory, the nonlinear Wiener structure provides strong model reduction capabilities beyond those of linear and bilinear models.
The Wiener structure has already been used implicitly in many works applying Koopman theory to autonomous systems, e.g., \cite{Takeishi.2017, Lusch.2018, Otto.2019}, without an explicit classification as such.
However, to our knowledge, Wiener-type Koopman models with controls have not yet been used in the literature.
Herein, we derive the Wiener-type Koopman form for input-affine systems, showing that the MIMO Wiener structure is naturally related to Koopman theory.

Our model identification and reduction procedure extends the Koopman deep-learning framework for autonomous systems by Lusch et al.~\cite{Lusch.2018} to MIMO models.
To this end, we combine autoencoders with non-autonomous linear dynamical blocks.
This extension is mathematically justified, as shown by our derivation of the Wiener-type Koopman form.
As opposed to existing MIMO Wiener identification methods
\cite{Westwick.1996, Lovera.2000, Arto.2001, 
Gomez.2005, Janczak.2007,
Hsu.2009, 
Schoukens.2011, Mu.2012}, 
our Koopman strategy targets data-driven model reduction, and therefore uses full-state information (trajectory snapshots) instead of input-output data.
The method is simple to use thanks to few requirements on model properties, identification procedure and data set, e.g., step responses are sufficient.

We apply our model reduction and identification framework to three case studies: a system with input multiplicity, a chemical reactor and a distillation column.
We demonstrate that the low-order MIMO Wiener Koopman models reproduce the dynamics of the nonlinear input-affine systems precisely.
Further, we compare to bilinear and linear models.
In particular, we show that low-order Wiener-type models reach similar or better accuracy than low-order bilinear models, while only involving a linear dynamical block.
Thus, Wiener-type models enable a higher degree of reduction and simplification.

The manuscript is structured as follows:
In \cref{sec:methods}, we provide the theoretical background on Koopman theory and block-structured modeling.
In \cref{sec:derivation}, we derive a Koopman model of Wiener type.
In \cref{sec:identification}, we present a deep-learning identification and model reduction strategy using Wiener-type Koopman or MIMO Wiener models.
In \cref{sec:casestudy}, we demonstrate our model identification framework in three case studies.
We conclude the manuscript with summary and outlook in \cref{sec:conclusion}.

\section{Theoretical Background}
\label{sec:methods}
We provide a brief review of system identification using block-structured modeling, and subsequently summarize Koopman theory and its extensions to control systems. 
Our presentation and notation of Koopman theory follows the comprehensive reviews
in \cite{Williams.2015, Mauroy.2020, Otto.2021,Brunton.2021}.

\subsection{Block-structured modeling}
Block-structured modeling is a data-driven system identification approach combining nonlinear static and linear dynamical blocks \cite{Giri.2010}. 
Block-structured models have been used extensively in control applications for several decades, and are prominent due to their intuitive structure and
sophisticated identification strategies available \cite{Giri.2010,Schoukens.2019}. 
Moreover, there exist efficient methods to embed serial block-structured models in MPC applications, e.g., \cite{Norquay.1998}. 

\begin{figure}[ht]
	\centering
	\begin{tikzpicture}[x=6ex, y=6ex, node distance=1ex and 5ex,
	squarednode/.style={rectangle, draw=black!100, fill=white!0, thin, minimum width=10ex, minimum height=6ex},]
	\node[squarednode, minimum width=3cm]      (dynamic)  {
	$\left.\ddt{\bm z}\right|_t  = A \bm z(t) + B \bm u(t)$
	};
	\node[font=\small] (text1) [below = of dynamic] {Linear Dynamics};
	\node[squarednode]      (static)       [right = of dynamic] {$\bm y(t) = \bm N(\bm z(t))$};
	\node[font=\small] (text1) [below = of static] {Static Nonlinearity};
	\draw[-latex] ($(dynamic.west) + (-1.,0.)$) -- (dynamic.west)
	node [pos=0.5,above,font=\small] {$\bm u(t)$};
	\draw[-latex] (dynamic.east) -- (static.west)
	node [pos=0.5,above,font=\small] {$\bm z(t)$};
	\draw[-latex] (static.east) -- ($(static.east) + (1.,0.)$) 
	node [pos=0.5,above,font=\small] {$\bm y(t)$};
	\end{tikzpicture}
	\caption{Block-structure of MIMO Wiener model. 
	Linear dynamics involve inputs $\bm{u}(t)\in \mathbb{R}^{n_u}$ and states $\bm{z}(t)\in \mathbb{R}^{n_z}$.
	Output nonlinearity $\bm{N}: \mathbb{R}^{n_z} \rightarrow \mathbb{R}^{n_y}$ yields outputs $\bm{y}(t) \in \mathbb{R}^{n_y}$.}
	\label{fig:wiener}
\end{figure}
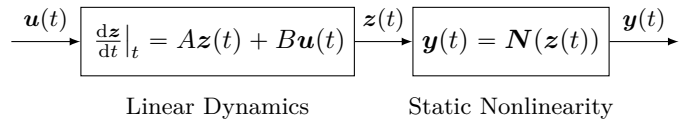

Wiener models are one common type of serial block-structured models, cascading a linear dynamical block with a nonlinear static block, Fig.~\ref{fig:wiener}.
Using MIMO blocks leads to MIMO Wiener models, 
which is regarded as an extension of the single-input single-output (SISO) parallel Wiener \cite{Schoukens.2012} structure. 
However, compared to SISO models, the mapping of MIMO blocks is mathematically more complex and usually associated with more parameters.
In many applications, Wiener structures are assumed to be purely empirical.
However, some approaches assign the model blocks a certain physical meaning (gray-box modeling),
e.g., the nonlinear static block representing the nonlinear stationary states \cite{Pearson.2000} or measurement nonlinearity. 
Furthermore, depending on the formulation, Wiener structures are regarded as a restriction
or a variant of Volterra series \cite{Boyd.1985}.
In \cref{sec:derivation}, we show that the Wiener structure can also be mathematically motivated by Koopman theory.

Within the Wiener model class, there exists a variety of options for formulating the nonlinearity and linear dynamical block. 
Besides using, e.g., polynomial or piecewise linear functions for the static nonlinear block, several works have investigated neural networks (NNs) to represent the nonlinear block in
SISO \cite{AlDuwaish.1996,
	Schram.1996, 
	Aadaleesan.2008, 
	Janczak.2003b, 
	Peng.2011} 
and MIMO \cite{Arto.2001, 
	Hsu.2009}
Wiener models. 
Combining discrete-time linear blocks and NNs promotes deep-learning environments as an alternative to standard model identification toolboxes.
In addition, identification algorithms for MIMO Wiener systems have been proposed, which can be
grouped into those relying on special input signals \cite{Westwick.1996,Lovera.2000, Schoukens.2011, Mu.2012}, in particular Gaussian excitation, and those presupposing invertibility of the static nonlinear output block \cite{Gomez.2005, Janczak.2007}.
Overall, MIMO Wiener identification methods rely on specific mathematical properties of the model blocks, expert knowledge or special characteristics of the identification data set.

\subsection{Koopman Theory}
Koopman theory facilitates the representation of nonlinear dynamics in a linear framework \cite{Mezic.2005, Mauroy.2020}.
In essence, Koopman theory applies nonlinear coordinate transformation to lift the system states into an infinite-dimensional space, where the dynamics behave linearly.
The Koopman operator acts on a function space $\mathcal{F}$,
e.g., Hilbert space, of nonlinear
observable functions $g \in \mathcal{F} : X \rightarrow \mathbb{C}$ of the states $\bm x(t)\in X \subseteq \mathbb{R}^{n_x}$.
Here, the Koopman operator advances these nonlinear observables in time in a linear fashion. 

Consider a nonlinear continuous-time autonomous system:
\begin{equation}
\label{eqn:nonlinear}
\ddt{\bm x}  = \bm f(\bm x) \,,
\end{equation}\noindent
where $\bm x(t)$ 
is the differential state at time $t$, 
and $\bm f: X \rightarrow \mathbb{R}^{n_x}$ 
are the continuously differentiable nonlinear
dynamics.
The generator $\lie_f : \mathcal{F} \rightarrow \mathcal{F}$, also referred to as Lie operator, belongs to the Koopman operator family and induces linear dynamics:
\begin{equation}
\label{eqn:conti_Koopman}
\ddt{g}  = \lie_f\, g \,,
\end{equation}\noindent
where: 
\begin{equation}
\lie_f\, g(\bm x) \equiv \nabla_x g(\bm x) \cdot \bm f(\bm x) \,.
\end{equation}\noindent
This means that the nonlinear dynamics become globally linear, but the lifted states, i.e., all possible nonlinear observables, are infinite dimensional.
Since the Koopman operator is linear, it has
eigenfunctions $\varphi_i: X \rightarrow \mathbb{C}$ with corresponding eigenvalues $\lambda_i \in \mathbb{C}$
evolving according to the characteristic eigendynamics:
\begin{equation}
\label{eqn:eigendynamics}
\ddt{\varphi_i}  = \lambda_i  
\varphi_i \,.
\end{equation}\noindent
If an observable $g_j$ can be represented in the coordinates spanned by these eigenfunctions,
its dynamic evolution is fully described by the eigendynamics. 

For the practical application of Koopman theory, 
a finite-dimensional invariant subspace $Z$ is sought, 
spanned by a finite subset of observables or eigenfunctions.
Combining a finite-dimensional 
matrix approximation $A: Z \rightarrow Z$ (projection) of the Lie operator on
this subspace $Z$ 
and the respective nonlinear coordinate lifting $\bm T: X \rightarrow Z$,
the continuous-time Koopman model with linear dynamics becomes:
\begin{subequations}
\label{eqn:koopmanmodel}
\begin{align}
\ddt{\bm z} &= A \bm z \,, \label{eqn:linear_dynamics}\\
\bm z(t) &= \bm T(\bm x(t)) \,,
\end{align}
\end{subequations}\noindent
in the observables $\bm z(t)\in Z$.
Similarly, the corresponding eigenfunctions $\bm \varphi(\bm x)$ satisfy:
\begin{equation}
\label{eqn:eigendynamics_invariant}
\ddt{\bm \varphi} = \Lambda \bm \varphi \,.
\end{equation}\noindent
We limit ourselves to systems with point spectrum \cite{Mezic.2020}
and assume $A$ and $\Lambda$ to be constant, thus linear time-invariant (LTI) dynamics.
In practical implementations, a ``Koopman canonical transform'' \cite{Surana.2016b} is often useful, where \cref{eqn:koopmanmodel,eqn:eigendynamics_invariant} are formulated in terms of real-valued observables and eigenfunctions, respectively.
We adopt this formulation hereinafter and assume that $\Lambda$
has block-diagonal form. 

For nonlinear systems exhibiting a single globally attracting  equilibrium point,  
a finite-dimensional Koopman invariant subspace can (theoretically) always be constructed \cite{Lan.2013}. 
Moreover, if the original states are included in the subspace $Z$, then $\bm x$ can be reconstructed from $\bm z$
by linear recombination, i.e., $\bm x = C \bm z$.
Conversely, for any system exhibiting multiple fixed points 
or more general attractors,
a 
finite-dimensional invariant subspace 
including the original states $\bm x$ cannot exist \cite{Brunton.2016b}.
In practice, 
determining 
a truly invariant low-dimensional subspace $Z$ with associated coordinate transformation $\bm T (\bm x)$ is impractical for most problems.
Instead, 
local approximations to invariant subspaces 
of the Koopman operator are constructed while trying to avoid 
closure issues. 
These approximations may also be interpreted as a truncated Koopman mode decomposition
\cite{Williams.2015}.
Often, they can provide sufficiently
accurate predictions of the nonlinear dynamics in short-term 
predictions on an application-relevant part of the state space.

The so-called dynamic mode decomposition (DMD) \cite{Schmid.2010, Tu.2013} and extended DMD \cite{Williams.2015} are two data-driven model identification methods closely related to Koopman theory.
While extended DMD essentially performs linear regression and modal decomposition on a dictionary of nonlinear observables $\bm z=\Theta(\bm x)$ generated from time-series data $\bm x(t)$, 
the original DMD uses the more restrictive identity dictionary $\bm z=\bm x$, 
trying to approximate \cref{eqn:koopmanmodel} or its discrete-time counterpart.

\subsection{Deep-learning strategies for Koopman models}
\label{sec:deeplearning}
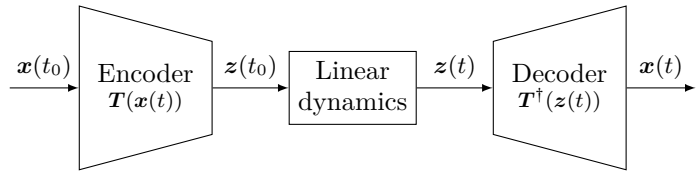
\begin{figure}[ht]
	\centering
	\scalebox{1.0}{
	\begin{tikzpicture}[x=6ex, y=6ex, node distance=1ex and 5ex,
	squarednode/.style={rectangle, draw=black!100, fill=white!0, thin, minimum width=10ex, minimum height=6ex},]
	\node[draw, trapezium, black, 
	trapezium left angle = 75, trapezium right angle = 75,
	rotate=-90, trapezium stretches body, 
	text width=1cm, align=center] at (0,0) (encoder) {\rotatebox{90}{\parbox[c]{10ex}{\centering Encoder \footnotesize$\bm T(\bm x(t))$}}};
	\node[squarednode, align=center, minimum width=10ex] at (3,0)      (dynamic) {Linear\\dynamics};
	\node[draw, trapezium, 
	trapezium left angle = 75, trapezium right angle = 75,
	black, rotate=90, trapezium stretches body, 
	text width=1cm, align=center] at (6,0) (decoder) {\rotatebox{-90}{\parbox[c]{10ex}{\centering Decoder \footnotesize$\bm T^{\dag}(\bm z(t))$}}};
	\draw[-latex] ($(encoder.south) + (-1.,0.)$) -- (encoder.south)
	node [pos=0.5,above,font=\small] {$\bm x(t_0)$};
	\draw[-latex] (encoder.north) -- (dynamic.west)
	node [pos=0.5,above,font=\small] {$\bm z(t_0)$};
	\draw[-latex] (dynamic.east) -- (decoder.north)
	node [pos=0.5,above,font=\small] {$\bm z(t)$};
	\draw[-latex] (decoder.south) -- ($(decoder.south) + (1.,0.)$) 
	node [pos=0.5,above,font=\small] {$\bm x(t)$};
	\end{tikzpicture}
	}
	\caption{Schematic illustration of deep-learning Koopman architectures for autonomous dynamical systems. Adapted from \cite{Lusch.2018}.
	}
	\label{fig:deeplearning}
\end{figure}

Several works have investigated deep learning to discover Koopman observables or EDMD dictionaries for finite-dimensional representations of autonomous dynamical systems
\cite{Li.2017, Takeishi.2017, Lusch.2018, Otto.2019, Yeung.2019, Erichson.2019, Pan.2020}.
Given the unknown form of the observables $\bm z(t)=\bm T(\bm x(t))$, 
using artificial neural networks is expedient 
due to their universal approximation capabilities \cite{Li.2017}.
While some approaches include the original states in the learned transformation 
\cite{Li.2017, Yeung.2019}, 
most works use autoencoder networks to allow for a nonlinear state reconstruction \cite{Takeishi.2017, Lusch.2018, Otto.2019, Erichson.2019, Pan.2020},
see \cref{fig:deeplearning}.
Moreover, restricting the dimension of the
learned Koopman subspace to a latent variable space,
combines Koopman theory and model reduction \cite{Lusch.2018, Otto.2019}.
In this respect, Koopman deep-learning architectures are related to other data-driven reduction strategies based on latent dynamics and nonlinear projection through autoencoders, e.g., \cite{Lee.2020, Gedon.2020, Tsay.2020b, Masti.2021}.
However, Koopman methods distinguish from 
these methods by simultaneously generating nonlinear transformation and linear latent dynamics.

A deep-learning framework as depicted in \cref{fig:deeplearning} readily provides approximations of the nonlinear transformation $\bm T^{}$ (encoder) along with the approximate inverse $\bm T^{\dag}$ (decoder).
While the autoencoders are typically trained on all data sampled from the relevant state space, the encoder mainly serves to provide transformed initial conditions for a simulation task, whereas the decoder recovers the original states along the dynamic trajectory.
Many works use uniformly sampled data sets and prefer to identify the mathematically equivalent discrete-time linear dynamics to \cref{eqn:linear_dynamics}:
\begin{equation}
    \bm z_{k+1} = \bar A \bm z_k,
\end{equation}\noindent
where $k=0,1,...$ are the sampling instants.
In the present literature, a simultaneous identification strategy for autoencoder and linear dynamics is achieved by summing
loss functions for the nonlinear reconstruction error as well as linear or nonlinear prediction errors over multiple time steps, e.g., \cite{Lusch.2018}.
Regularization strategies, such as Tikhonov or $\ell_1$ regularization, are frequently used in the training to prevent 
overfitting and stability issues,
while closure issues of the subspace can be detected 
through validation.
In addition, sparsity-promoting training for a block-diagonal structure of the linear dynamics or expert knowledge for presetting the structure of $A$ can foster the specific discovery of the eigendynamics, \cref{eqn:eigendynamics_invariant}, and supports dynamical decoupling and uniqueness of the encoding and decoding \cite{Lusch.2018, Pan.2020}. 
For example, Lusch \etal\cite{Lusch.2018} combine deep autoencoders with a linear block-diagonal state-space model as depicted in \cref{fig:deeplearning}.
\subsection{Koopman theory for controlled systems}
\label{sec:Koopman_affine}
Classical Koopman theory was proposed for autonomous systems and is not directly extendable to causal nonlinear systems with external forcing of the general form:
\begin{equation}
\label{eqn:control}
\ddt{\bm x} = \bm f'(\bm x, \bm u) \,.
\end{equation}\noindent
Herein, $\bm f' : X \times U \rightarrow \mathbb{R}^{n_x}$ additionally depends on the inputs $\bm u(t) \in U \subseteq \mathbb{R}^{n_u}$.

Empirical extension of the linear Koopman representation of an autonomous system, \cref{eqn:koopmanmodel}, by external inputs can 
corrupt the linear dynamics \cite{Brunton.2016b, Otto.2021}. 
While 
efforts have been made to cast representations of dynamical systems with inputs into a Koopman-operable form \cite{Williams.2016, Proctor.2018, Korda.2018, Peitz.2019, Klus.2020, Kaiser.2021},
no canonical form of \cref{eqn:control} with Koopman dynamics explicit in the external inputs has been presented in the literature yet.
However, Surana \cite{Surana.2016} and Goswami \& Paley \cite{Goswami.2017} showed that
a universal Koopman form exists for the continuous-time input-affine nonlinear systems class:
\begin{equation}
\label{eqn:control_affine}
\ddt{\bm x}  = \bm f(\bm x) + \sum_{i=1}^{n_u} \bm h_i(\bm x) u_i \,,
\end{equation}\noindent
where $\bm f(\bm x)$ and $\bm h_i(\bm x)$ act on the state space.
This class represents many physical systems and engineering problems.
In addition, more general systems, \cref{eqn:control}, can be brought into this form when embedding additional input dynamics for $\bm u(t)$.
Since our derivation of a Wiener-type Koopman model in \cref{sec:derivation} builds on the works of Surana \cite{Surana.2016}, we discuss their derivation and follow-up works in detail below.

Assuming that all states $\bm x$ can be linearly reconstructed from a finite number of Koopman observables or eigenfunctions, i.e., $\bm x = C^{\varphi}\bm \varphi(\bm x)$,
Surana \cite{Surana.2016} demonstrated that the Koopman dynamics for systems of the form (\ref{eqn:control_affine}) are also input-affine:
\begin{equation}
\label{eqn:koopmanmodel_controls}
\ddt{\bm\varphi}  = \Lambda \bm \varphi + \sum_{i=1}^{n_u} \nabla_x \bm \varphi^T \, \bm h_i(\bm x) \, u_i \, \rvert_{\bm x = C^{\varphi}\bm \varphi} \,. \\
\end{equation}\noindent
where, noticeably, $u_i$ are still the original inputs rather than transformed ones.
Here, the eigenfunctions spanning the Koopman-invariant subspace take the form $\bm \varphi^\mathrm{c}(\bm x,\bm u)=\bm \varphi(\bm x)$, which are eigenfunctions associated with the autonomous dynamical part, $\bm f(\bm x)$.
\Cref{eqn:koopmanmodel_controls} underlines that a prerequisite for obtaining a Koopman-invariant subspace, 
i.e., an exact finite-order representation of the dynamics,
is to find eigenfunction coordinates $\bm \varphi(\bm x)$ such that $\nabla_x \bm \varphi(\bm x)^T \bm h_i(\bm x)$ lie in the span of $\bm \varphi(\bm x)$ \cite{Goswami.2017, Surana.2016}:
\begin{equation}
\label{eqn:Surana}
\nabla_x \bm \varphi(\bm x)^T \bm h_i(\bm x) \equiv \sum_{j=1}^{N} \varphi_j(\bm x)\, \bm v^{(h_i)}_j\,,
\end{equation}\noindent
where $\bm v^{(h_i)}_j$ are the vector-valued Koopman modes from a decomposition associated with $\bm h_i(\bm x)$.
In practice, this assumption 
can be hard to satisfy, and determining a closed representation may require a significantly higher number of Koopman eigenfunctions
than needed for the autonomous counterpart, e.g., \cite{Otto.2021}.
Finally, \cref{eqn:koopmanmodel_controls} 
can be brought into a form useful for control \cite{Surana.2016, Goswami.2017}: 
\begin{equation}
\label{eqn:koopman_eigendynamics_bilinear}
\ddt{\bm\varphi} = \Lambda \bm \varphi 
+ \sum_{i=1}^{n_u} \Xi^{(i)} \bm \varphi u_i   \,,
\end{equation}\noindent
where the last term is bilinear in $u_i$ and $\varphi_j$, and $\Lambda$ and $\Xi^{(i)}$ are matrices.
Alternatively, the more intuitive Koopman model in terms of the nonlinear observables $\bm z$ reads: 
\begin{equation}
\label{eqn:koopmanmodel_controls_bilin_Surana}
\begin{split}
\ddt{\bm z} &= A\bm z + \sum_{i=1}^{n_u} B^{(i)} \bm z u_i  \,, \\
\bm x &= C \bm z \,,\\ 
\bm z_0 &= \bm T(\bm x_0) \,.
\end{split}
\end{equation}\noindent
Herein, for the sake of clarity, we explicitly state the linear reconstruction, the initial condition, where $\bm x_0 = \bm x(0)$, and denote the nonlinear coordinate lifting by $\bm T(\bm x(t))$.
This Koopman form employs a bilinear dynamical model and linear state reconstruction.

Surana \cite{Surana.2016} notices that the Koopman decomposition may require infinitely many observables for an arbitrarily accurate representation of the input-affine dynamics.
Indeed, finite-order bilinear and linear models are known to be incapable of exhibiting chaotic behavior or output multiplicity \cite{Pearson.2003,Brunton.2016b}.
While the general limits of linear and bilinear models are known \cite{Pearson.2003}, a rigorous classification of systems that can be arbitrarily exactly approximated by means of a finite-order bilinear Koopman form is missing. 
However, in many cases, a finite-order Koopman model may represent a truncated or approximate decomposition.
Note further that the bilinearity of the Koopman model, \cref{eqn:koopmanmodel_controls_bilin_Surana}, only holds for the continuous-time generator $\lie_f$ but not generally in discrete time \cite{Peitz.2020}.
While an exact transformation between continuous-time and discrete-time bilinear models is possible for (pairwise) separable systems \cite{Dunoyer.1997,Rumschinski.2012},
a discrete-time bilinear model may otherwise only be used for sufficiently fine sampling or weak nonlinearity.

According to \cref{eqn:koopmanmodel_controls}, the bilinear dynamics reduce to linear dynamics, if there exist eigencoordinates satisfying \cite{Otto.2021}:
\begin{equation}
\label{condition:linear}
\nabla_x \bm \varphi(\bm x)^T h_i(\bm x) \equiv \mathrm{const.}
\end{equation}\noindent
If Condition (\ref{condition:linear}) holds, we obtain a Koopman representation
with linear dynamics: 
\begin{equation}
\label{eqn:koopmanmodel_controls_lin_Surana}
\begin{split}
\ddt{\bm z} &= A \bm z + B \bm u \,, \\
\bm x &= C \bm z \,,\\
\bm z_0 &= \bm T(\bm x_0) \,.
\end{split}
\end{equation}\noindent
Note that \cref{eqn:koopmanmodel_controls_lin_Surana} is closely related to the DMD with controls \cite{Proctor.2016} approach.

Korda and Mezi\'{c} \cite{Korda.2018} proposed to empirically restrict the Koopman model to a linear form, \cref{eqn:koopmanmodel_controls_lin_Surana}.
Several works employed a linear Koopman approach in identification and control studies and observed a higher performance over conventional linear identification methods
\cite{Korda.2018, Arbabi.2018, Narasingam.2019, Abraham.2019, Folkestad.2020, Han.2020, Son.2021b, Ping.2021}.
However, when significant nonlinearity is present, \cref{eqn:koopmanmodel_controls_lin_Surana} was found to lead to prediction errors and tracking losses
\cite{Peitz.2020,Bruder.2021,Folkestad.2021}.
The works \cite{Peitz.2020,Bruder.2021,Folkestad.2021} demonstrated the higher predictive ability of bilinear models, \cref{eqn:koopmanmodel_controls_bilin_Surana}, 
a result known from standard system identification literature \cite{Pearson.2003}.
Deep learning for Koopman models with inputs was applied in \cite{Abraham.2019, Han.2020, Ping.2021}.
However, these works identified linear models, \cref{eqn:koopmanmodel_controls_lin_Surana}, without a nonlinear decoder.
\section{Koopman control model with invertible coordinate transformation}
\label{sec:derivation}
In this section, we derive a Koopman model with Wiener-type block structure.
In the derivation of the global bilinearization, \cref{eqn:koopmanmodel_controls_bilin_Surana}, Surana \cite{Surana.2016} assumes that $\bm \varphi(\bm x)$ is injective and that $\bm  x$ can be fully and uniquely recovered by linear recombination of Koopman eigenfunctions.
However, this assumption can become restrictive when aiming at low-dimensional, possibly linear, latent dynamical models \cite{Surana.2016, Kaiser.2021}.
The identification of underlying dominant nonlinear patterns, as targeted by nonlinear model reduction, may be difficult or impossible when requiring linear state reconstruction, cf.~also\cite{Lee.2020}.
On the other hand, due to the exact convertibility of LTI dynamics with zeroth-order hold between discrete-time and continuous-time as well as the rich control and system theory, linear dynamics are practically favorable over bilinear dynamics.

We aim to unify the simplicity of linear dynamics, the accuracy of bilinear Koopman representations, and the model-reduction capabilities by nonlinear state reconstruction.
To this end, we drop the linear state recovery assumption of \cite{Surana.2016} and assume instead that the coordinate transformations $\bm \varphi (\bm x)$ and $\bm T(\bm x)$ are continuously invertible, i.e., continuously bijective.
Instead of \cref{eqn:koopmanmodel_controls}, we formulate:
\begin{equation}
\label{eqn:koopmanmodel_controls_invertible}
\ddt{\bm\varphi}  = \Lambda \bm \varphi + \sum_{i=1}^{n_u} \nabla_x \bm \varphi^T \, \bm h_i(\bm x) \, u_i \, 
\rvert_{\bm x = \bm \varphi^{-1}(\bm \varphi)} \,.
\end{equation}\noindent
Furthermore, the requirement described by \cref{eqn:Surana} still applies.
Analogue to \cref{eqn:koopmanmodel_controls_bilin_Surana}, we then obtain a new Koopman model with nonlinear state reconstruction:
\begin{equation}
\label{eqn:koopmanmodel_controls_bilin}
\begin{split}
\ddt{\bm z} &= A\bm z + \sum_{i=1}^{n_u} B^{(i)} \bm z u_i  \,, \\
\bm x &= \bm T^{-1}(\bm z) \,,\\
\bm z_0 &= \bm T(\bm x_0) \,.
\end{split}
\end{equation}\noindent
While \cref{eqn:koopmanmodel_controls_bilin_Surana} and \cref{eqn:koopmanmodel_controls_bilin} appear similar, the matrices and coordinate transformations may be fundamentally different.
Under Condition (\ref{condition:linear}), this model reduces to the linear-dynamics form:
\begin{equation}
\label{eqn:koopmanmodel_controls_lin}
\begin{split}
\ddt{\bm z} &= A \bm z + B \bm u \,, \\
\bm x &= \bm T^{-1}(\bm z) \,,\\
\bm z_0 &= \bm T(\bm x_0) \,.
\end{split}
\end{equation}\noindent
In the following, we will use this model form with latent linear dynamics.
The Koopman model \cref{eqn:koopmanmodel_controls_lin} is a block-structured MIMO Wiener model, cf.~\cref{fig:wiener}.
Notice that this model form extends to the MIMO Hammerstein-Wiener structure when adding an input nonlinearity block, i.e., when enriching the Koopman observables by nonlinear functions of the inputs.

The mathematical derivation of the Koopman model, \cref{eqn:koopmanmodel_controls_lin}, also supports the block-structure of the MIMO Wiener and parallel Wiener classes and provides 
a mathematical interpretation of these model architectures.
On the other hand, the known limitations of Wiener modeling \cite{Boyd.1985, Pearson.2003} can be applied to analyze the model.
In particular, Boyd and Chua \cite{Boyd.1985} find that the Wiener-type structure 
is able to capture the behavior of systems with ``fading memory'' to arbitrary accuracy.
In fading memory systems, the current state only depends on recent input signals but not on the distant past, which is related to a unique steady state \cite{Boyd.1985}.
In practice, the fading memory property may be easier to judge than the existence of a finite-dimensional subspace satisfying Condition (\ref{condition:linear}).

\section{Identification Strategy}
\label{sec:identification}
\begin{figure}[h!]
	\centering
	\scalebox{1.0}{
	\begin{tikzpicture}[x=6ex, y=6ex, node distance=1ex and 5ex,
	squarednode/.style={rectangle, draw=black!100, fill=white!0, thin, minimum width=10ex, minimum height=6ex},]
	\node[draw, dashed, trapezium, 
	trapezium left angle = 75, trapezium right angle = 75,
	black, rotate=-90, trapezium stretches body, 
	text width=1cm, align=center] at (0,0) (encoder) {\rotatebox{90}{\parbox[c]{10ex}{\centering Encoder \footnotesize$\bm T(\bm x_j)$}}};
	\node[squarednode, align=center, minimum width=10ex] at (3,0)      (dynamic) {Linear\\dynamics};
	\node[draw, trapezium, 
	trapezium left angle = 75, trapezium right angle = 75,
	black, rotate=90, trapezium stretches body, 
	text width=1.cm, align=center] at (6,0) (decoder) {\rotatebox{-90}{\parbox[c]{10ex}{\centering Decoder \footnotesize$\bm T^{\dag}(\bm z_{k})$}}};
	\draw[-latex, dashed] ($(encoder.south) + (-1.,0.)$) -- (encoder.south)
	node [pos=0.5,above,font=\small] {$\bm x_j$};
	\draw[-latex] (encoder.north) -- (dynamic.west)
	node [pos=0.5,above,font=\small] {$\bm z_j$};
	\draw[-latex] (dynamic.east) -- (decoder.north)
	node [pos=0.5,above,font=\small] {$\bm z_k$};
	\draw[-latex] (decoder.south) -- ($(decoder.south) + (1.,0.)$) 
	node [pos=0.5,above,font=\small] {$\bm x_k$};
	\draw[-latex] ($(dynamic.south) - (0.,0.6)$) -- (dynamic.south) 
	node [pos=0.4,right,font=\small] {$\bm u_k$};
	\end{tikzpicture}
	}
	\caption{Deep-learning strategy for MIMO Wiener-type models. Original states $\bm x$, transformed states $\bm z$, time instants $j,k = 0,1,...$, and $k\geq j$.}
	\label{fig:encoder_controls_linear}
\end{figure}
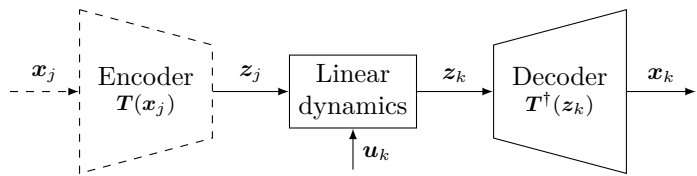

In this section, we present a deep-learning identification procedure for the Wiener-type Koopman model type with external inputs that we derived in the previous section.
The proposed method can also be regarded as a MIMO Wiener identification strategy.
Our method is motivated by the Koopman deep-learning framework for autonomous systems by Lusch et al.~\cite{Lusch.2018}.
However, the proposed procedure differs from the work of \cite{Lusch.2018} in that we consider systems with external inputs, \cref{fig:encoder_controls_linear}.
The derivation given in \cref{sec:derivation} provides a mathematical foundation for this extension.

\Cref{fig:encoder_controls_linear} illustrates our NN block-structured model, which implements a discrete-time dynamical block sandwiched by a deep autoencoder, cf.~\cref{eqn:koopmanmodel_controls_lin}.
We use the autoencoder to construct the nonlinear state transformation, $\bm T(\,\cdot\,)$, along with its approximate inverse, $\bm T^{\dag}(\,\cdot\,)$.
The encoder is used in the identification procedure but is not necessarily part of a control model. 
The number and arrangement of neurons as well as the order $n_z$ are model hyperparameters and affect the accuracy and degree of reduction.
Similar to \cite{Lusch.2018}, we choose to identify discrete-time LTI dynamics because they do not need a computation of continuous time-derivatives in the training, and can be directly embedded as algebraic equations in an optimization problem.
If desired, the transformation into continuous-time representation may be conducted after the training.
Our procedure can be easily implemented within a deep-learning environment, e.g., Tensorflow or PyTorch.
The identification does not rely on special input signals and enables a simultaneous identification of all model blocks.
Importantly, although the nonlinear state transformation blocks are invertible,
the strategy can be used to identify systems with input multiplicity in individual outputs as long as the steady state is unique, cf.~\cref{sec:casestudy1}.

To train the model, we use snapshots of step-response trajectories of all states with sufficiently fine sampling.
Since we target data-driven model reduction, we assume that noise-free full-state information is available.
We subdivide this data into shorter trajectories containing $p$ snapshots and randomly group multiple trajectories in training batches.
The trajectory length $p$ and the batch size are training hyperparameters.
Therefore, the provided identification data can be a collection of short samples.
Finally, we split the batched data set into training and validation data and use an independent test data set.
In the training, we encode the initial state snapshot $\bm x_0$ of a trajectory and decode single and multi-step predictions generated by the dynamical block.
To this end, we sum loss terms that we compute for every $p$-samples trajectory in a training batch using the mean squared error (MSE):
\begin{subequations}
\begin{align}
	\label{eqn:loss1}
L_1 =&  \frac{1}{p} \sum_{k=0}^{p} \lVert \bm x_k - \bm T^{\dag}(\bm T(\bm x_k))\rVert\ind{MSE} \;,\\
L_2 =&  \frac{1}{p-1} \sum_{k=0}^{p-1} \lVert \bm x_{k+1} - \bm T^{\dag}(\bm z_{k+1}(\bm x_k)) \rVert\ind{MSE} \label{eqn:loss2} \;,\\
L_3 =&  \frac{1}{p-1} \sum_{k=0}^{p-1} \lVert \bm x_{k+1} - \bm T^{\dag}(\bm z_{k+1}(\bm x_0)) \rVert\ind{MSE} \label{eqn:loss3} \;,
\end{align}
\end{subequations}\noindent
where single and multi-step predictions $\bm z_{k+1}(\bm x_j)$ are computed as:
\begin{equation}
\label{eqn:loss_pred}
\begin{split}
\bm z_{k+1} &= A \bm z_{k} + B \bm u_{k} \,,\hspace{10pt} k=j,j+1,...\\
\bm z_j &= \bm T(\bm x_j) \,.\\
\end{split}
\end{equation}\noindent
We aim at discovering the eigendynamics and thus restrict the matrix $A\in\mathbb{R}^{n_z \times n_z}$ to a block-diagonal structure.
The term $L_1$ evaluates the nonlinear reconstruction accuracy of the autoencoder. 
The second term, $L_2$, considers single-step predictions starting from various points along each trajectory.
The third term, $L_3$, evaluates multi-step predictions of up to $p$ time steps to account for slow dynamics, all starting from the initial state $\bm x_0$ of a trajectory.
We scale all losses by the mean squared sum of the data points, yielding normalized mean squared error (NMSE) loss.
Additionally, we add $\ell_1$ regularization on all trainable parameters.
We train on a weighted linear combination of all loss terms.
Details on the implementation are provided with the case studies.

\section{Case Studies}
\label{sec:casestudy}
We demonstrate our identification strategy in three case studies with identification problems from the literature.
The first case study investigates a nonlinear system with input multiplicity \cite{Brunton.2016b}.
Next, we consider the chemical reactor investigated by \cite{Narasingam.2019}.
Finally, the third case study performs model identification of a distillation column \cite{Pearson.2000}.

In all case studies, we train Wiener-type Koopman models using our identification strategy from \cref{sec:identification}.
We implement the identification framework with Python 3.9 using the Tensorflow 2.5 \cite{Tensorflow} API.
We provide all code of our identification framework online\footnote{See: \url{https://git.rwth-aachen.de/avt-svt/public/koopman-wiener}}.
To identify suitable hyperparameters, 
we perform an a-priori systematic parameter search by training a variety of model structures including shallow and deep autoencoder structures with varying order $n_z$ of the latent linear
dynamics.
For the sake of simplicity, we use a symmetric autoencoder design for the block-structured Koopman model, i.e., the arrangement and number of neurons in the hidden layers of the decoder is reversed as in the encoder.
Further details on the training procedure and hyperparameters are provided in the Appendix.

We identify linear as well as bilinear dynamical models in all case studies as a benchmark and use the Koopman deep-learning framework described here.
To this end, we remove the nonlinear hidden layers of the decoder while keeping the nonlinear encoding, resulting in a linear recombination of the latent coordinates.
These benchmark models are discrete-time versions of \cref{eqn:koopmanmodel_controls_bilin_Surana} and \cref{eqn:koopmanmodel_controls_lin_Surana}. 
The bilinear dynamics include an explict term linear in the controls, as this was found to be beneficial in the training. 
We use sufficiently fine sampling to ensure the validity of discrete-time bilinear dynamics \cite{Peitz.2020}.

In all case studies, the training data are generated from numerical simulations.
To use single-platform Python code throughout this work, all dynamic models are implemented in Pyomo.DAE \cite{PyomoDAE} and simulated with CasADi \cite{Casadi}.
We assume noise-free snapshot data.
We evaluate the identified models on an independent test data set comprising dynamic step-response trajectories that were not seen in the training.
In the test, we perform a complete forward simulation given only the initial state $\bm x_0$ and a sequence of inputs $\bm u_k$ over the full simulation horizon.
We compare the predictions of the identified models to the true output profiles from the mechanistic model, and judge the performance by means of the NMSE over all scaled model outputs at all sampling instants.
As a benchmark, we also identify MIMO Wiener models using the MATLAB 2021a System Identification Toolbox \cite{matlab} using the same model architectures as described in the case studies, default settings and trained on the scaled data set. 
Since the MATLAB toolbox failed to generate low-order models of acceptable accuracy,
we show these results only in the Appendix.

\subsection{Case Study I: System with input multiplicity}
\label{sec:casestudy1}
We investigate a variant of the input-affine nonlinear system studied by Brunton \etal\cite{Brunton.2016b}:
\begin{equation}
\label{eqn:casestudy1}
\begin{split}
\ddt{x_1} &= -0.1 x_1 + u\,,\\
\ddt{x_2} &= x_1^2 - x_2\,,\\
y &= x_2 \,.
\end{split}
\end{equation}\noindent
The system is stable and exhibits input multiplicity, since $y\ind{s} = 100 u^2$.
However, it has a unique steady state, $\bm x\ind{s} = [10u, 100u^2]^T$.
As shown in \cite{Brunton.2016b}, an analytic bilinear Koopman model, \cref{eqn:koopmanmodel_controls_bilin_Surana}, can be constructed for this example when specifying $\bm z := [x_1, x_2, x_1^2]^T$.
Here, we demonstrate that block-structured Wiener-type Koopman models with even fewer differential states $\bm z$ can be identified using our deep-learning framework. 

\subsubsection{Data sampling and training}
The data used for the identification is generated in a dynamical simulation and comprises the trajectories of all states resulting from a sequence of 100 random uniform input steps with $u\in[-1,1]$.
Each input step has a duration of $t\ind{step} = 200$, and the time grid is sampled with $\Delta t = 1$.
In the same way, we create a test data set comprising 20 random uniform input steps and $t\ind{step} = 100$.
We identify full-state models, i.e., $\bm y = \bm x$.
The hyperparameter search indicated that low-order latent dynamics, $n_z =2$, can reproduce the dynamics sufficiently accurately.
Moreover, shallow encoders/decoders with a single hidden layer of $20$ neurons are a valid choice.

\begin{figure}[htb]
\centering
	\begin{subfigure}[]{0.48\textwidth}
		\centering
		\includegraphics[width=\textwidth]{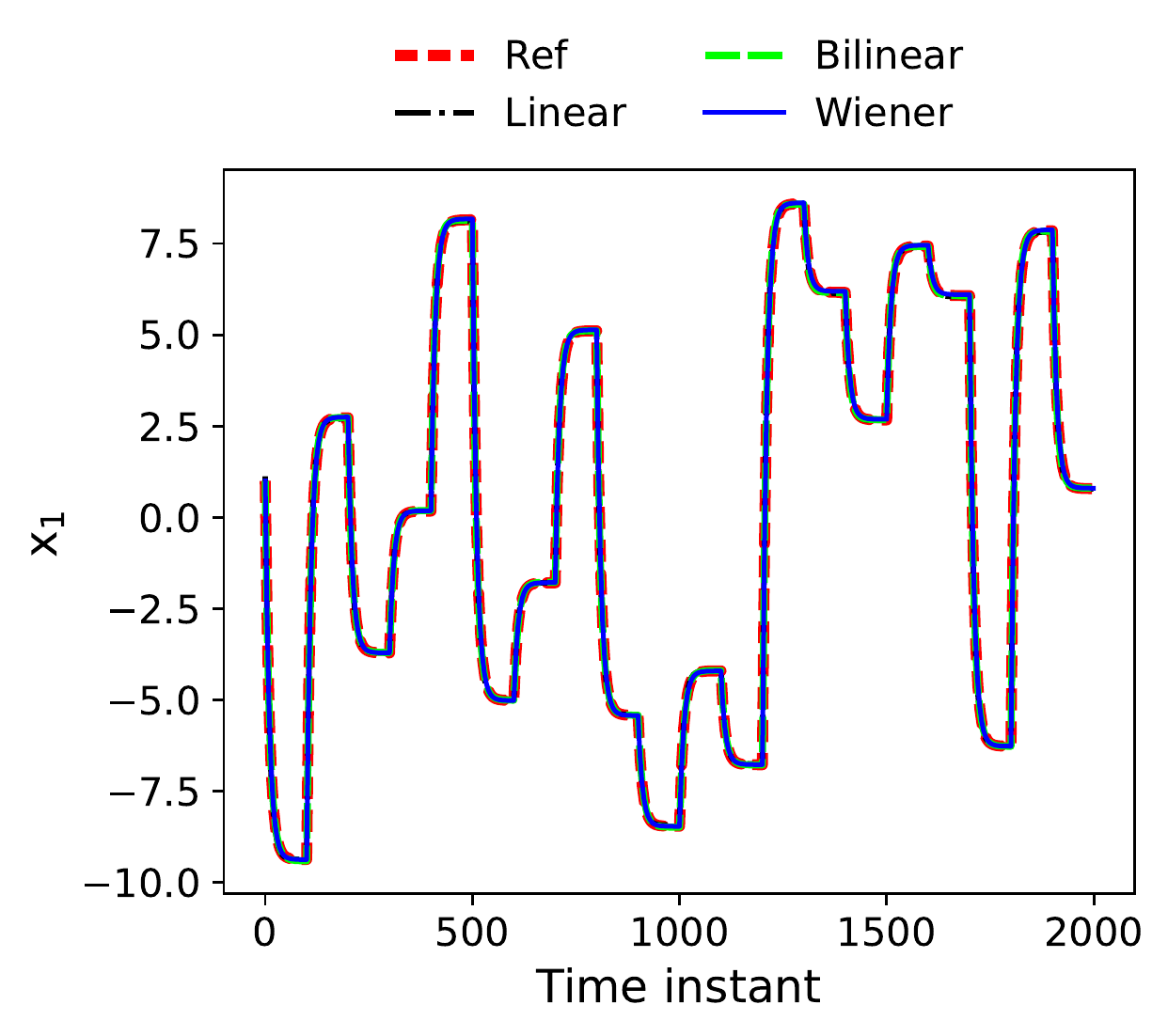}
		\caption{}
		\label{fig:casestudy_1_nz2a}
	\end{subfigure}
	\begin{subfigure}[]{0.45\textwidth}
		\centering
		\includegraphics[width=\textwidth]{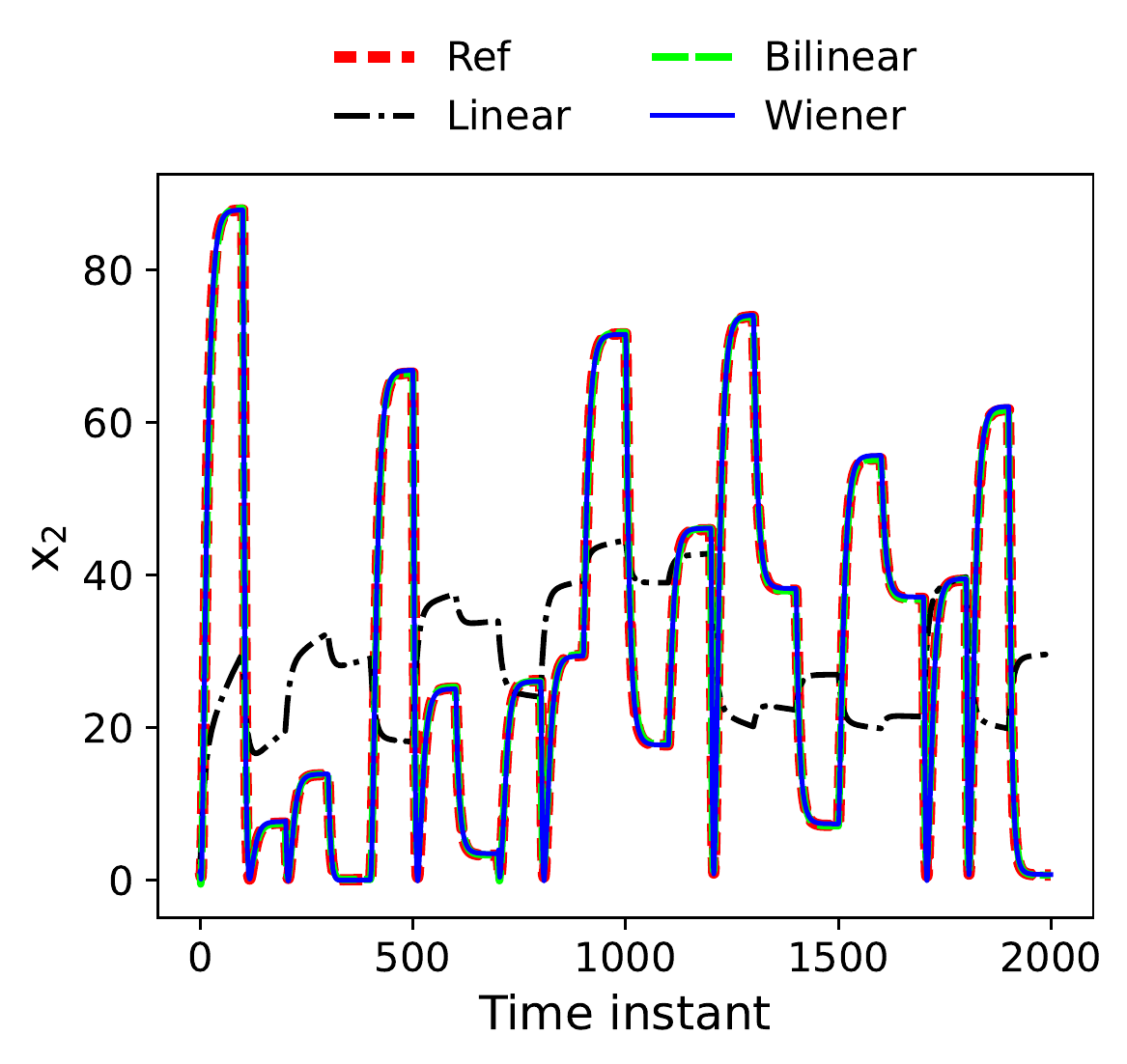}
		\caption{}
		\label{fig:casestudy_1_nz2b}
	\end{subfigure}
	\begin{subfigure}[]{0.45\textwidth}
		\centering
		\includegraphics[width=\textwidth]{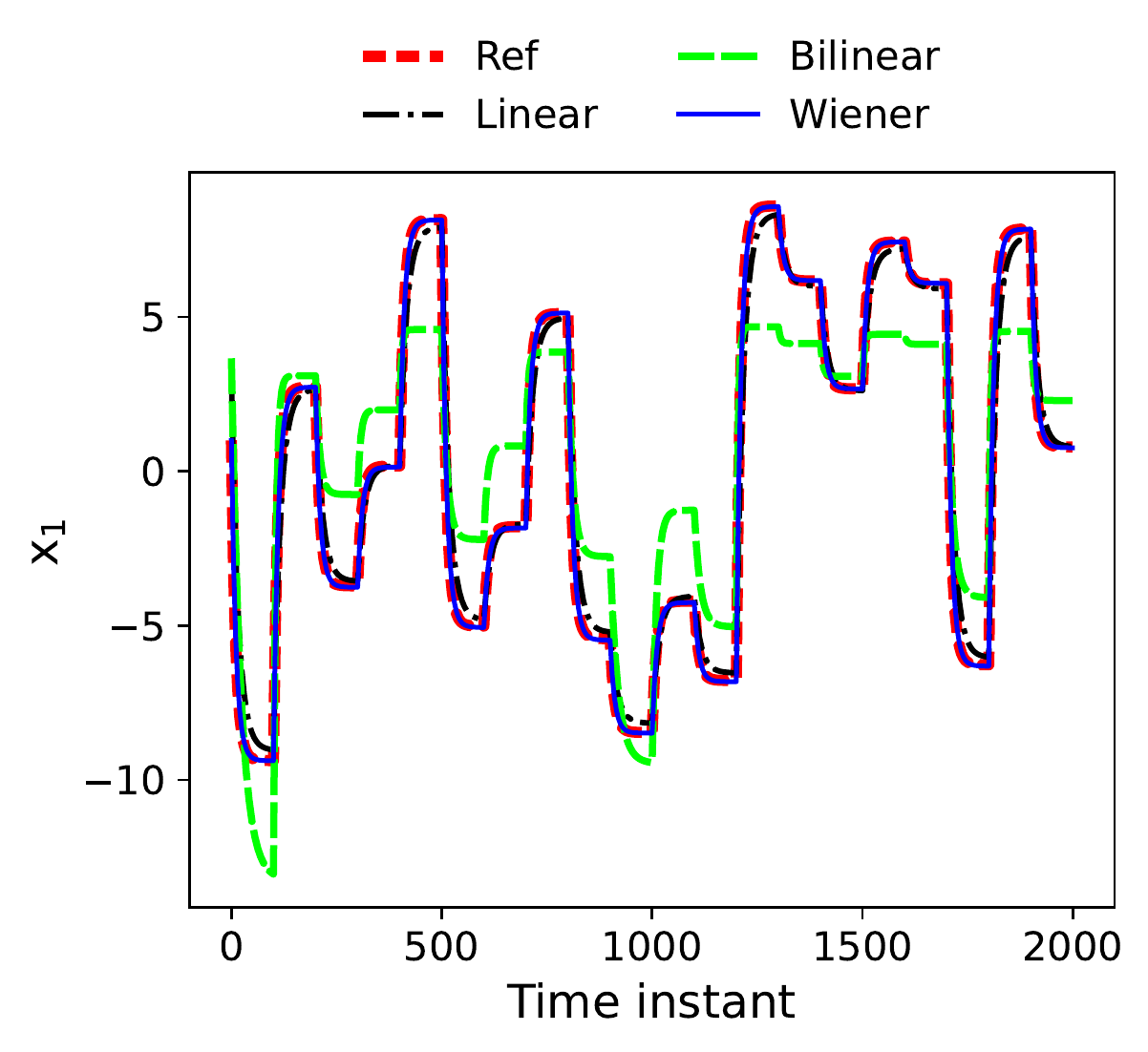}
		\caption{}
		\label{fig:casestudy_1_nz1a}
	\end{subfigure}
	\begin{subfigure}[]{0.45\textwidth}
		\centering
		\includegraphics[width=\textwidth]{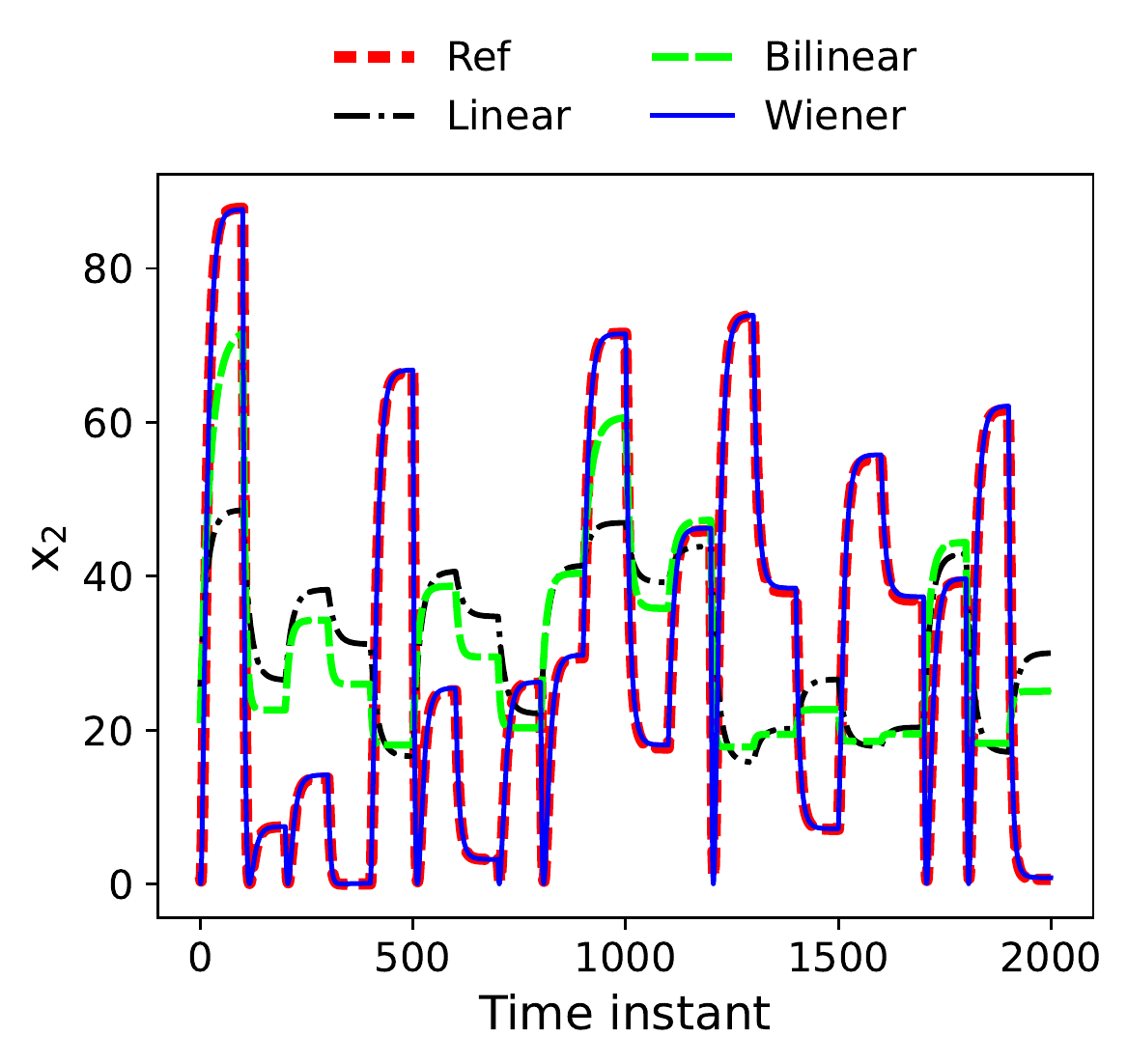}
		\caption{}
		\label{fig:casestudy_1_nz1b}
	\end{subfigure}
	\caption{Case study 1: Performance of the identified full-state models with $n_z=2$ (a, b) and $n_z=1$ (c, d) on the test set. 
	}
	\label{fig:results_casestudy1}
\end{figure}
\subsubsection{Results and Discussion}
\Cref{fig:casestudy_1_nz2a,fig:casestudy_1_nz2b} compare the block-structured Koopman model, linear and bilinear models for $n_z=2$.
\Cref{tab:hyperparameters_1} summarizes the model hyperparameters and prediction errors for all models.
All models reproduce the $x_1$ dynamics precisely, \cref{fig:casestudy_1_nz2a}.
Furthermore, both the Wiener-type and bilinear model provide an accurate prediction of $x_2$, whereas the linear model fails at reproducing the respective dynamics, \Cref{fig:casestudy_1_nz2b}.
As also evident from \cref{tab:hyperparameters_1}, increasing
the linear order to $n_z=10$ degrades the predictions rather than improving the model.
\begin{table}[ht]
	\caption{\small Case study 1: Hyperparameters and prediction errors of the MIMO models: Subspace dimension $n_z$, no.~of hidden neurons of encoder $n\ind{n,E}$ and decoder $n\ind{n,D}$,
	no.~of original states $n_x = 2$.}
	\label{tab:hyperparameters_1}
	\centering
	\begin{tabular}{l l l l l }
		\toprule
		\textbf{Model} & $\bm{n_z}$ & $\bm{n\ind{n,E}}$ & $\bm{n\ind{n,D}}$ 
		&  \textbf{NMSE} \\ 
		\midrule
		Wiener  & 2 & 20 & 20 
		& 2$\cdot$10$^{-5}$\\[1ex]
		Bilinear & 2 & 20 & $-$ 
		& 2$\cdot$10$^{-4}$\\[1ex]
		\multirow{2}{*}{Linear} 
		&  2 & 20 & $-$ 
		& 0.15\\[0ex]
		& 10 & 50 & $-$ 
		& 0.77 \\
		\bottomrule
	\end{tabular}
\end{table}

Since we target model reduction, we further reduce the latent space to $n_z=1$ while keeping the autoencoder complexity, and investigate the performance of the identified models, \cref{fig:casestudy_1_nz1a,fig:casestudy_1_nz1b}.
Now, both linear and bilinear model fail to capture the system response, while the Wiener-type model continues to predict the dynamic and stationary behavior successfully.
This indicates that the system dynamics allow for a joint encoding-decoding of both states, i.e., nonlinear projection, while linear reconstruction is not applicable.
These results underline the strong model reduction capabilities of the Wiener structure.

\subsection{Case Study II: Chemical reactor}
In this case study, we investigate the continuously operated chemical reactor for the exothermic reaction A $\rightarrow$ B presented in \cite{Alanqar.2015, Narasingam.2019}:
\begin{equation}
\label{eqn:casestudy2}
\begin{split}
\ddt{c\ind{A}} &= \frac{F}{V} (c\ind{A,in} - c\ind{A}) - k_0 \exp\left(-\frac{E}{RT}\right) c\ind{A}^2
\,,\\
\ddt{T} &= \frac{F}{V} (T\ind{in} - T) - \frac{k_0 \Delta H}{\rho c_p} \exp\left(-\frac{E}{RT}\right) c\ind{A}^2 + \frac{\dot{Q}}{\rho c_p V}
\,.
\end{split}
\end{equation}\noindent
Herein, the molar concentration $c\ind{A}$ and the reactor temperature $T$ are the differential states $\bm x$.
Similar to \cite{Narasingam.2019}, we consider the heat duty $\dot{Q}$ as the varying system input $u$, and assume all other parameters as time-invariant.
Further details on the model and parametric values can be found in \cite{Alanqar.2015, Narasingam.2019}. 
Most importantly, the system exhibits three steady states for $\dot{Q} = 0$. 
Narasingam \& Kwon \cite{Narasingam.2019} identified a linear Koopman model with $n_z = 7$ for the high-yield operating region of product B, and presented a stabilizing MPC to control the process.
However, the authors required expert knowledge to specify a dictionary of Koopman observables rather than employing a universal learning strategy.
While linear, bilinear and Wiener-type models are known to be incapable of representing systems with steady-state multiplicity globally exactly \cite{Pearson.2003}, we aim at a local model for the high-yield operating region.

\subsubsection{Data sampling and training}
We perform a dynamical simulation subject to a random sequence of 500 steps uniformly distributed over
$ \dot{Q} \in [-2000,10\,000] \,\unit{kJ\,h^{-1}}$ at a sampling rate of $\unit[1]{min^{-1}}$.
Each input steps has a duration of $t\ind{step} = \unit[4]{h}$.
Similarly, we create a test data set comprising 15 random uniform input steps with \unit[2]{h} duration.
We identify a model describing both states.
A systematic parameter screening suggested that accurate predictions can be attained with extremely parsimonious models using low-dimensional static and dynamic blocks.
Since we aim for low-order models rather than the most accurate predictor, we set $n_z = 1$ and use a shallow encoder/decoder each with a single hidden layer with 20 neurons.

\begin{figure}[htb]
	\begin{subfigure}[]{0.45\textwidth}
		\centering
		\includegraphics[width=\textwidth]{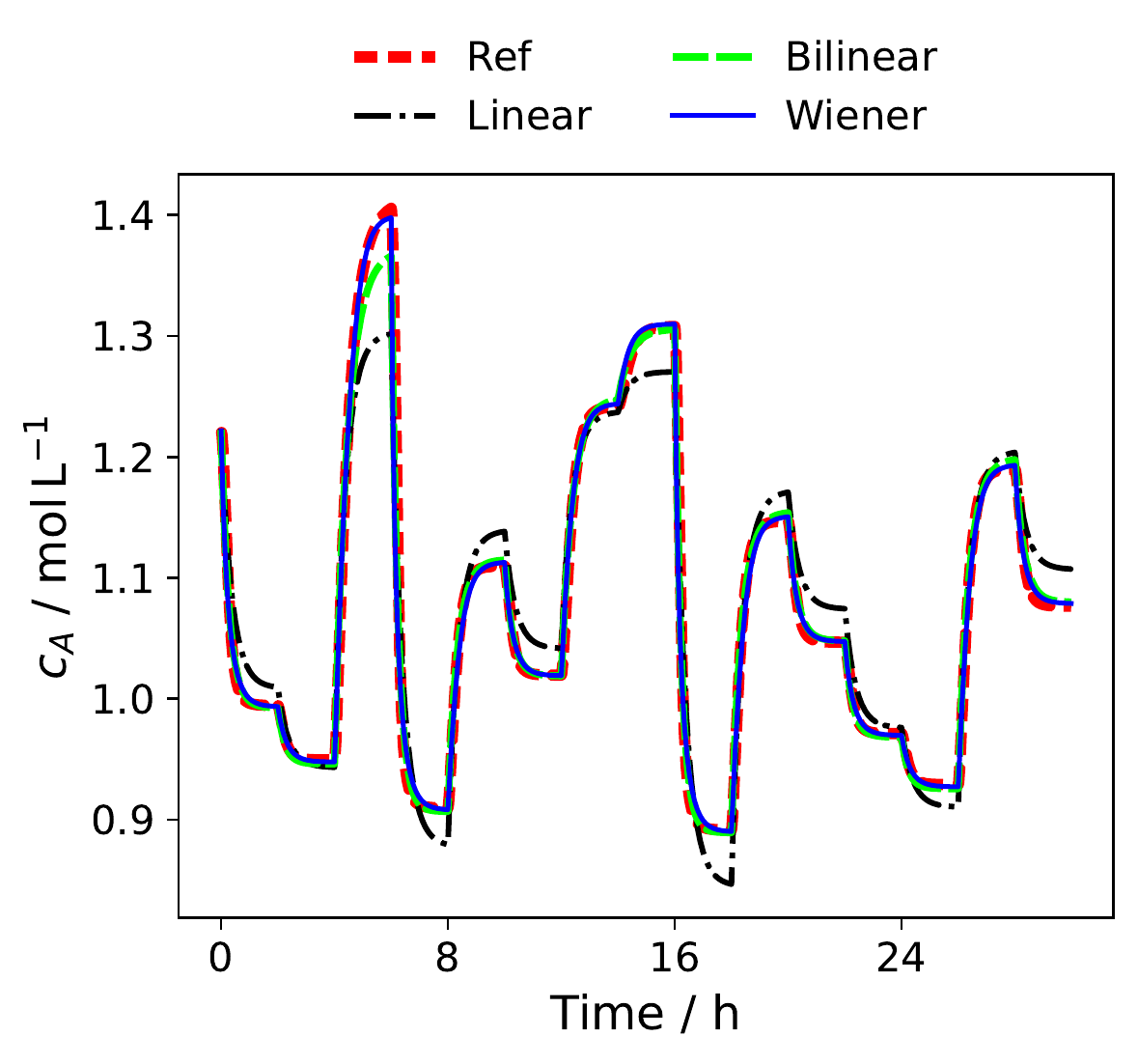}
		\caption{}
	\end{subfigure}
	\hfill
	\begin{subfigure}[]{0.45\textwidth}
		\centering
		\includegraphics[width=\textwidth]{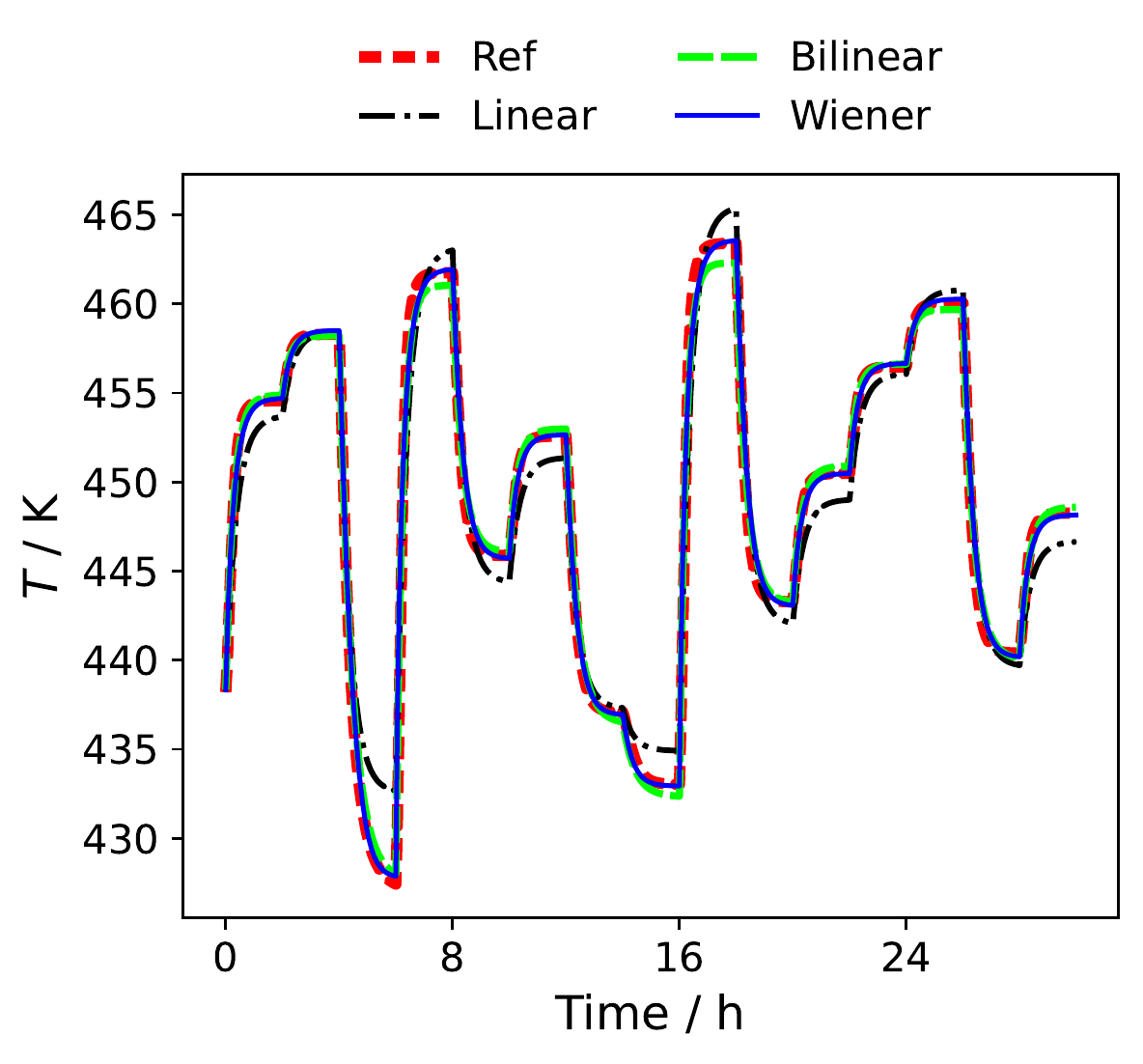}
		\caption{}
	\end{subfigure}
	\caption{Case study 2: Performance of the identified models on the test set. (a) Concentration of the reactant, (b) reactor temperature.}
	\label{fig:results_casestudy2}
\end{figure}
\subsubsection{Results and Discussion}
\Cref{fig:results_casestudy2} displays the model performance on the test set.
\Cref{tab:hyperparameters_2} collects the model parameters and the NMSE on the test set.
We find that the Koopman deep-learning framework enables the identification of low-order models of all model types.
However, despite the relatively narrow operating region, the prediction error of the linear model is noticeable and deviations occur as both over- and underestimation, \cref{fig:results_casestudy2}.
Besides steady-state offset of the linear model, there exists a pronounced error in the dynamical response.
Increasing the order of the linear model to $n_z=10$ slightly reduces the mismatch, but the prediction error is still present, see \cref{tab:hyperparameters_2}.
Especially, the steady-state offset does not vanish when increasing $n_z$.
In contrast, both Wiener and bilinear model successfully generate a precise forecast over the entire \unit[80]{h} horizon.
While the low-order bilinear model provides a slightly better dynamical response, the Wiener model is the only model with negligible steady-state offset over the entire operating range.
We conclude that by the simple means of introducing the shallow NN decoder block with few neurons for nonlinear reconstruction, a data-driven model with high accuracy is obtained.
Noticeably, although its linear dynamical block is significantly more simple than the original nonlinear dynamics, the low-order Wiener-type model reproduces the dynamics almost exactly.
The model could thus be used for highly precise tracking of quality requirements, allowing for smaller safety margins than a linear model.
\begin{table}[ht]
	\caption{\small Case study 2: Hyperparameters and prediction errors of the MIMO models: Subspace dimension $n_z$, no.~of hidden neurons of encoder $n\ind{n,E}$ and decoder $n\ind{n,D}$,
	no.~of original states $n_x = 2$.}
	\label{tab:hyperparameters_2}
	\centering
	\begin{tabular}{l l l l l}
		\toprule
		\textbf{Model} & $\bm{n_z}$ & $\bm{n\ind{n,E}}$ & $\bm{n\ind{n,D}}$ 
		&  \textbf{NMSE} \\ 
		\midrule
		Wiener  & 1 & 20 & 20 
		& 0.002\\[1ex]
		Bilinear   & 1 & 20 & $-$ 
		& 0.002\\[1ex]
		\multirow{2}{*}{Linear} 
		&  1 & 20 & $-$ 
		& 0.009\\[0ex]
		& 10 & 50 & $-$ 
		& 0.007 \\
		\bottomrule
	\end{tabular}
\end{table}
	
\subsection{Case Study III: Distillation Column}
\label{sec:casestudy_column}
High-purity distillation columns are one of the standard unit operations in chemical engineering.
The transient concentration profiles inside a distillation column exhibit a nonlinear traveling wave character, i.e., an underlying coherent structure moving inside the column \cite{Marquardt.1990}.
Block-structured SISO Wiener modeling of distillation columns has been the subject of several studies 
\cite{Schram.1996, Zhu.1999, Pearson.2000, Bloemen.2001, Norquay.1999}.
Further, the identification of linear models using log-scaled compositions is common practice in control engineering \cite{Skogestad.1988}, which is regarded as another variant of Wiener modeling.
In contrast, we aim at the identification of a low-order MIMO control model describing the evolution of all column states subject to multiple inputs.

In this case study, we investigate the high-purity methanol-propanol distillation column studied by Pearson \& Pottmann \cite{Pearson.2000}.
The column has $N=8$ equilibrium trays plus a total condenser and reboiler at the column top and bottom, respectively.
The binary feed mixture enters the column on the fourth stage and is purified in the rectifying section above and the stripping section below.
Details on the physical model and system parameters can be found in \cite{Pearson.2000}.
We employ our deep-learning framework for model reduction by identifying low-order MIMO dynamical models from sampled step-response state trajectories.
The two input variables are the feed composition of the heavy-boiling component ($u_1 = x_f$)
and the liquid molar reflux from the condenser into the column ($u_2 = L$).
Similar to \cite{Pearson.2000}, we assume a constant vapor flow rate $V$, constant feed flow rate $F$, and perfect control of the condenser holdup through the distillate flow rate $D$ in the model identification.
Hence, we do not consider these variables as model inputs.
The model has input-affine structure and $n_x=10$ dynamical states.

\subsubsection{Data sampling and training}
We consider the high-purity operating region of the column, characterized by small molar fractions, $1-x_{10}$, of the heavy-boiling component in the distillate stream.
The data used for the identification is generated in a dynamical simulation, and comprises the trajectories of all state responses to a random sequence of 400 input steps.
The respective input combinations $(u_1,u_2)$ lie on a uniform grid within the bounded set $x_f\in[0.5, 0.6]$ and $L\in[0.0155,0.0175]\,\unit{kmol/min}$.
Each input step has \unit[2]{h} duration at a data sampling rate of $\unit[1]{min^{-1}}$, providing a balanced amount of transient and steady-state data.
The simulation is initialized at the steady state corresponding to $\bm u = [0.6, 0.0175]^T$.

Since the molar fractions on the column trays vary over different orders of magnitude, 
we adopt the common practice \cite{Pearson.2000} to train on log-transformed molar fractions of the light-boiler in the stripping section and heavy-boiler in the rectifying section.
Due to the known nonlinear wave propagation phenomenon, we expect that a very low-order model is capable of providing satisfactory predictions.
This was confirmed by the parameter screening.
Hence, we identify MIMO models with a differential order $n_z=2$ and autoencoder networks with a single hidden layer, \cref{tab:hyperparameters}.
\begin{table}[ht]
	\caption{\small Case study 3: Hyperparameters and prediction error of the MIMO models: Subspace dimension $n_z$, no.~of hidden neurons of encoder $n\ind{n,E}$ and decoder $n\ind{n,D}$,
	no.~of original states $n_x = 10$.}
	\label{tab:hyperparameters}
	\centering
	\begin{tabular}{l l l l l}
		\toprule
		\textbf{Model} & $\bm{n_z}$ & $\bm{n\ind{n,E}}$ & $\bm{n\ind{n,D}}$ 
		&  \textbf{NMSE} \\ 
		\midrule
		Wiener & 2 & 10 & 10 
		& 0.0003\\[1ex]
		Bilinear  & 2 & 20 & $-$ 
		& 0.0012\\[1ex]
		\multirow{2}{*}{Linear} 
		 & 2 & 20 & $-$ 
		 & 0.0025\\[0ex]
		 & 10 & 50 & $-$ 
		 & 0.0013 \\
		\bottomrule
	\end{tabular}
\end{table}

\subsubsection{Results and Discussion}
\begin{figure}[htb]
	\centering
	\begin{subfigure}[]{0.45\textwidth}
		\centering
		\includegraphics[width=\textwidth]{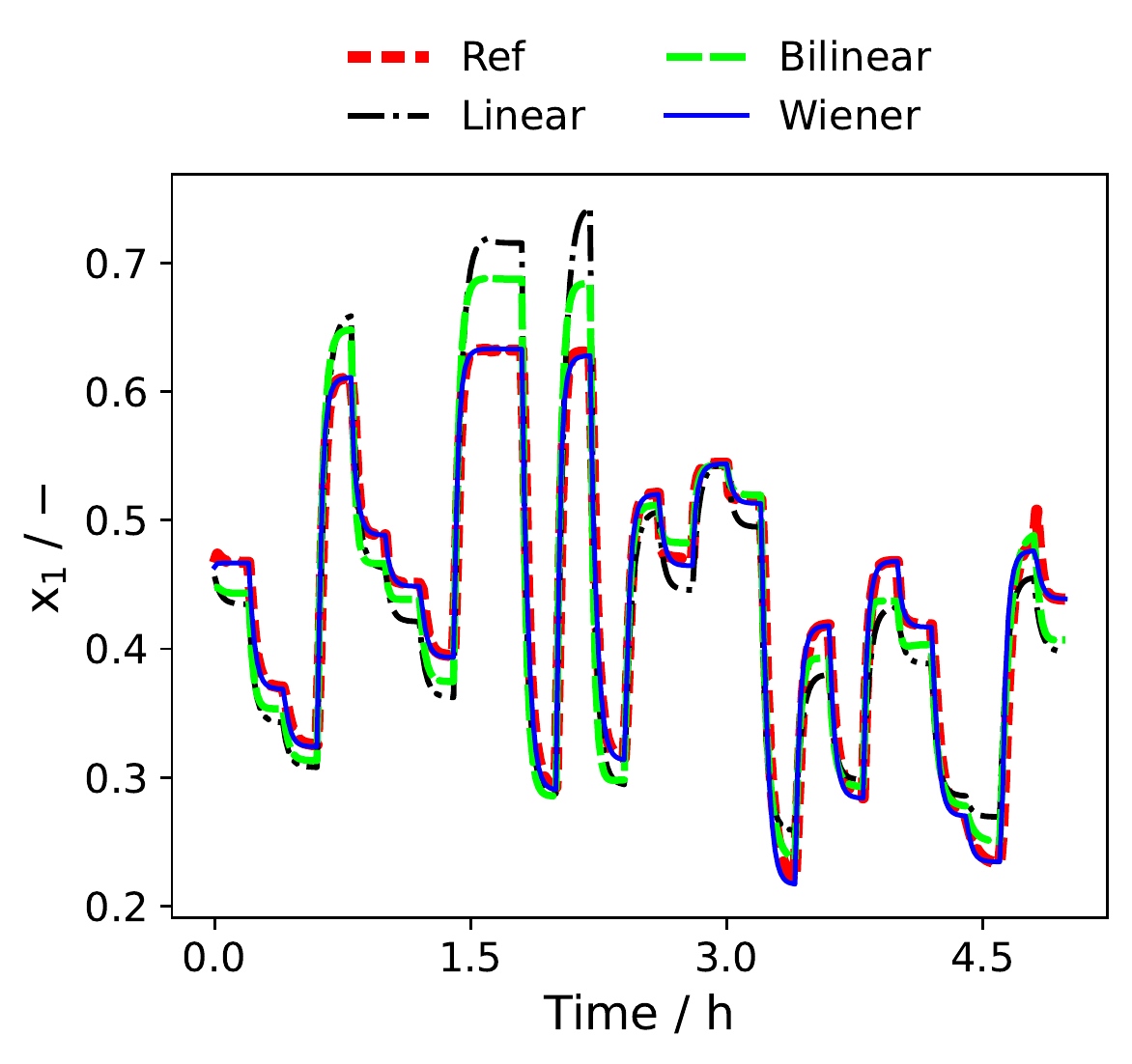}
		\caption{}
		\label{fig:bottomcomposition}
	\end{subfigure}
	\hfill
	\begin{subfigure}[]{0.48\textwidth}
		\centering
		\includegraphics[width=\textwidth]{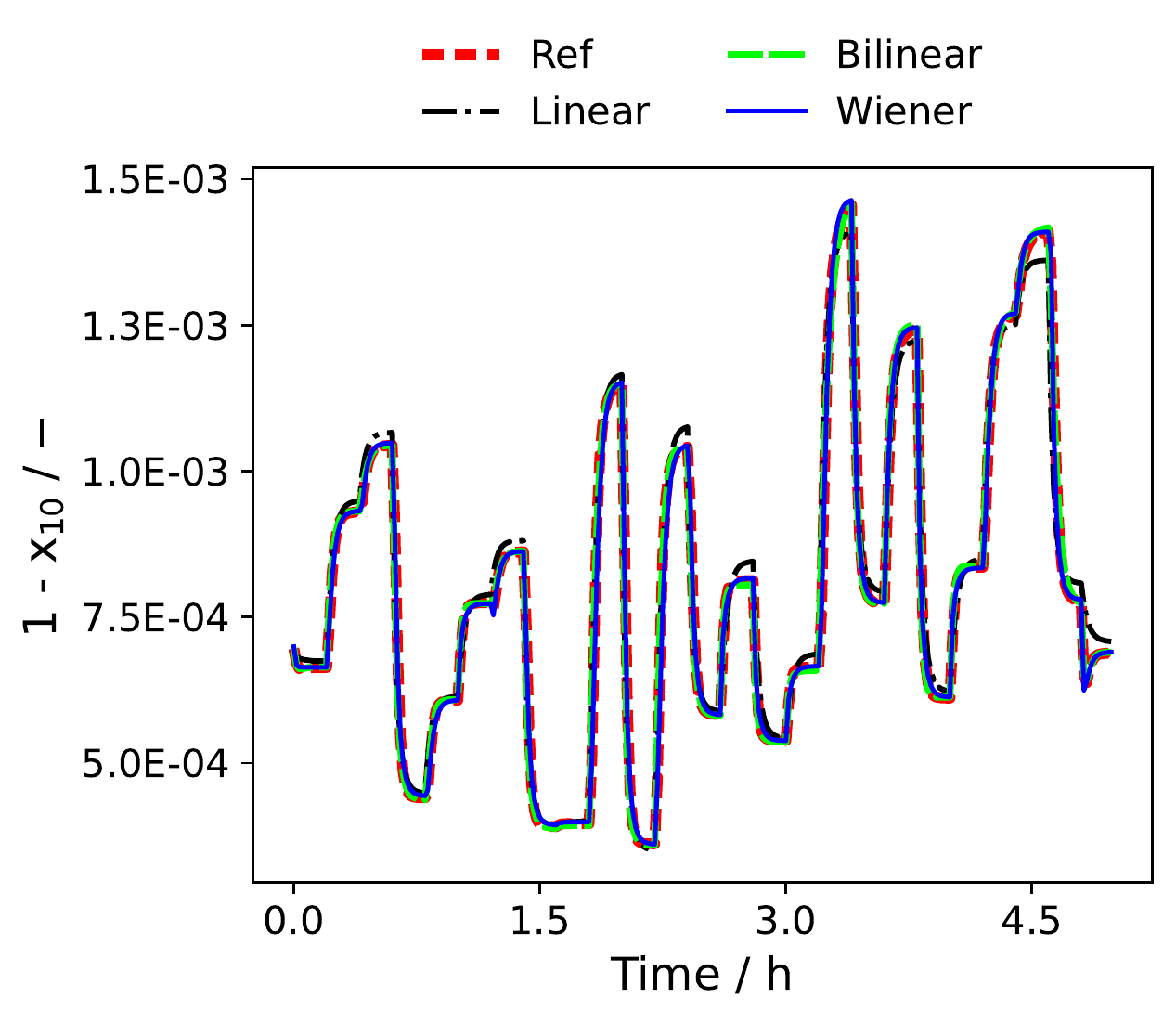}
		\caption{}
		\label{fig:topcomposition}
	\end{subfigure}
	\caption{Case study 3: Performance of the identified MIMO models on the test set. Predictions of the composition in the reboiler (a), and in the condenser (b).}
	\label{fig:results}
\end{figure}

\Cref{tab:hyperparameters} collects the values of the model hyperparameters along with the NMSE on the test set.
Moreover, \cref{fig:results} visualizes the predictions by the model structures.
The figure depicts the composition at the column bottom and top, as they characterize the exiting streams and are regarded as most critical with respect to prediction errors.
Moreover, the prediction errors throughout the column are of the same magnitude.
As visible from the figure, the Koopman deep-learning strategy enables the identification of accurate low-order control models.
All models, including the linear Koopman model trained on log-transformed data, are able to capture the dominant trends of the distillation column, and especially of the distillate composition.
However, the MIMO Wiener model still outperforms the linear and bilinear models.
In particular, \cref{fig:bottomcomposition} displays transient and steady-state offset in the predictions of the bottom composition by the linear and bilinear models.
Notice that increasing the complexity of the linear model results in a small improvement in the prediction quality but also a significantly higher number of trainable parameters, \cref{tab:hyperparameters}.
Moreover, notice that bilinear models provided increasingly accurate predictions when increasing the dynamic order, i.e., model complexity.
Conversely, even the low-order MIMO Wiener model yields highly precise long-term predictions for all column trays, which is also reflected in the NMSE performance index. 

\newpage
\section{Conclusions and Outlook}
\label{sec:conclusion}
We have presented a data-driven nonlinear model reduction and identification strategy for low-order modeling for control.
To this end, we have employed a class of MIMO Wiener-type Koopman models with external inputs.
While other works have focused on linear and bilinear Koopman models for input-affine systems, we have derived the Wiener-type Koopman variant here.
Our derivation also demonstrates that the MIMO Wiener model structure is naturally related to Koopman theory.
The Wiener structure is particularly suitable as a general-purpose model type for data-driven non-intrusive reduction, 
as it allows to identify dominant nonlinear dynamical patterns.

To obtain low-order MIMO control models, we have used a conceptually simple deep-learning strategy.
Therein, we pair autoencoders with linear dynamics, resulting in a Wiener-type Koopman model to learn the nonlinear dynamical behavior.
By restricting the latent space, i.e., linear dynamic block, to a low dimension, the deep-learning procedure is able to automatically generate strongly reduced models.
Further, employing a machine learning approach reduces the required expert knowledge significantly compared to a manual selection of Koopman observables or the application of intrusive nonlinear model reduction methods.
In contrast to other MIMO Wiener identification strategies in the literature, our framework uses state trajectory snapshots, implements a simultaneous identification of all model blocks and does not rely on special excitation signals.

We have discussed the advantages of the Wiener type over linear and bilinear Koopman models.
The Wiener block-structure unifies the simplicity of linear dynamical models and the accuracy of bilinear forms, and provides strong model reduction capabilities.
Due to the nonlinear reconstruction block, Wiener-type models can facilitate a more extensive reduction at higher accuracy than low-order bilinear models, although (high-order) bilinear models constitute a more universal model type.
Moreover, employing an LTI block in a Wiener structure enables
the application of standard methods for system analysis, allows for the exact transformation between discrete-time and continuous-time dynamics, and can be favorable in control and optimization procedures.
While a substitution of the linear dynamical block of the Wiener structure by more universal descriptors, e.g., \cite{Gedon.2020, Masti.2021}, is generally possible, such extensions lead to a higher number of trainable parameters and their numerical treatment and analysis is more difficult.
Moreover, our results suggest that using a MIMO Wiener model with linear dynamical block is already sufficient in many applications and yields a very minimalistic structure.

In three case studies, we have successfully demonstrated the identification and reduction capabilities of our deep-learning framework.
Throughout all case studies, 
low-order Wiener-type models outperformed the linear and bilinear Koopman models in terms of accuracy and feasible degree of reduction.
In addition, a standard identification toolbox failed to identify MIMO Wiener models for the systems and degree of reduction investigated here.
These results demonstrate that, in contrast to the common practice of identifying multiple MISO models, a suitable identification and reduction method can provide a single and extremely low-order MIMO model describing all states.

We conclude that data-driven model reduction by means of applied Koopman theory provides an effective and accurate framework for low-order modeling of moderately nonlinear control systems.
Therein, the MIMO Wiener structure constitutes a well suited and broadly applicable generic model form.
Future work should investigate Koopman-based reduction frameworks for strongly nonlinear, e.g., chaotic, systems with controls.
Moreover, our autoencoder framework could be 
extended to Koopman identification from noisy measurement data \cite{Mauroy.2016} and the identification in cases of missing state information \cite{Korda.2018, Kamb.2020}.

\begin{acknowledgements}
We gratefully acknowledge the  financial support of the Kopernikus project SynErgie 2 by the Federal Ministry of Education and Research (BMBF) and the project supervision by the project management organization Projekttr\"ager J\"ulich, as well as from the Helmholtz Association of German Research Centers as part of the Helmholtz School for Data Science in Life, Earth and Energy (HDS-LEE).
The authors thank Johannes M.~M.~Faust and Pascal Sch\"afer for proofreading and valuable feedback.
\end{acknowledgements}

\bibliographystyle{ieeetr}

\appendix

\section*{Appendix}
\subsection*{Additional information on the training}
We formulate the training problem according to \cref{sec:identification}.
Since we do not expect oscillatory open-loop behavior of the systems investigated, we restrict the system matrix $A$, \cref{eqn:loss_pred}, to a purely diagonal structure.
We use the ELU activation function and a linear output layer for the autoencoder networks.
All weights are initialized employing variance scaling with untruncated normal distribution, and all biases with zeros.
We scale the training data from zero to one individually for each differential state.
We specify the weights of the loss terms such that the three losses $L_1$, $L_2$ and $L_3$ are of the same order of magnitude in the later epochs.
The models are trained using batches of 32 trajectories, and 10,000 epochs to prevent premature termination.
We randomly divide the training data set into \unit[80]{\%} training data, \unit[20]{\%} validation data.
We execute all computations on a desktop computer with CPU at \unit[3.0]{GHz}.
We do not pre-train the encoder as we find no beneficial effect on the training result.
We employ the Adam \cite{Adam} training algorithm.
The model parameters with the smallest validation loss are saved and used afterwards.
Subsequently, we restore the model parameters form the epoch with the smallest validation loss.
To reduce stochastic effects, we train all structures multiple
times using different seeds.
\Cref{tab:training} collects the training hyperparameters.
Training on trajectories of length $s=50$ was found to yield good results, whereas shorter trajectories occasionally caused offset in both transients and stationary points.

\begin{table}[ht]
	\small
	\caption{Hyperparameters of the training problem.}
	\label{tab:training}
	\centering
	\begin{tabular}{p{5.cm} ll}
		\toprule
		\textbf{Parameter} & \textbf{Symbol} & \textbf{Value} \\ 
		\midrule
		Weight reconstruction & $\omega_1$ & 0.1  \\[1ex]
		Weight single-step prediction & $\omega_2$ & 1.0  \\[1ex]
		Weight multi-step prediction & $\omega_3$ & 1.0  \\[1ex]
		Weight $\ell_1$ regularization & $\omega\ind{r}$ & 1e-9  \\[1ex]
		Learning rate $ $ &&0.001 \\[1ex]
		Trajectory length & $p$ &50 \\[1ex]
		Batch size && 32 \\[1ex]
		No.~epochs && 10\,000 \\
		\bottomrule
	\end{tabular}\\
\end{table}

\newpage
\subsection*{Comparison to MATLAB System Identification Toolbox}
\begin{figure}[htb]
	\centering
	\begin{subfigure}[]{0.46\textwidth}
		\centering
		\includegraphics[width=\textwidth]{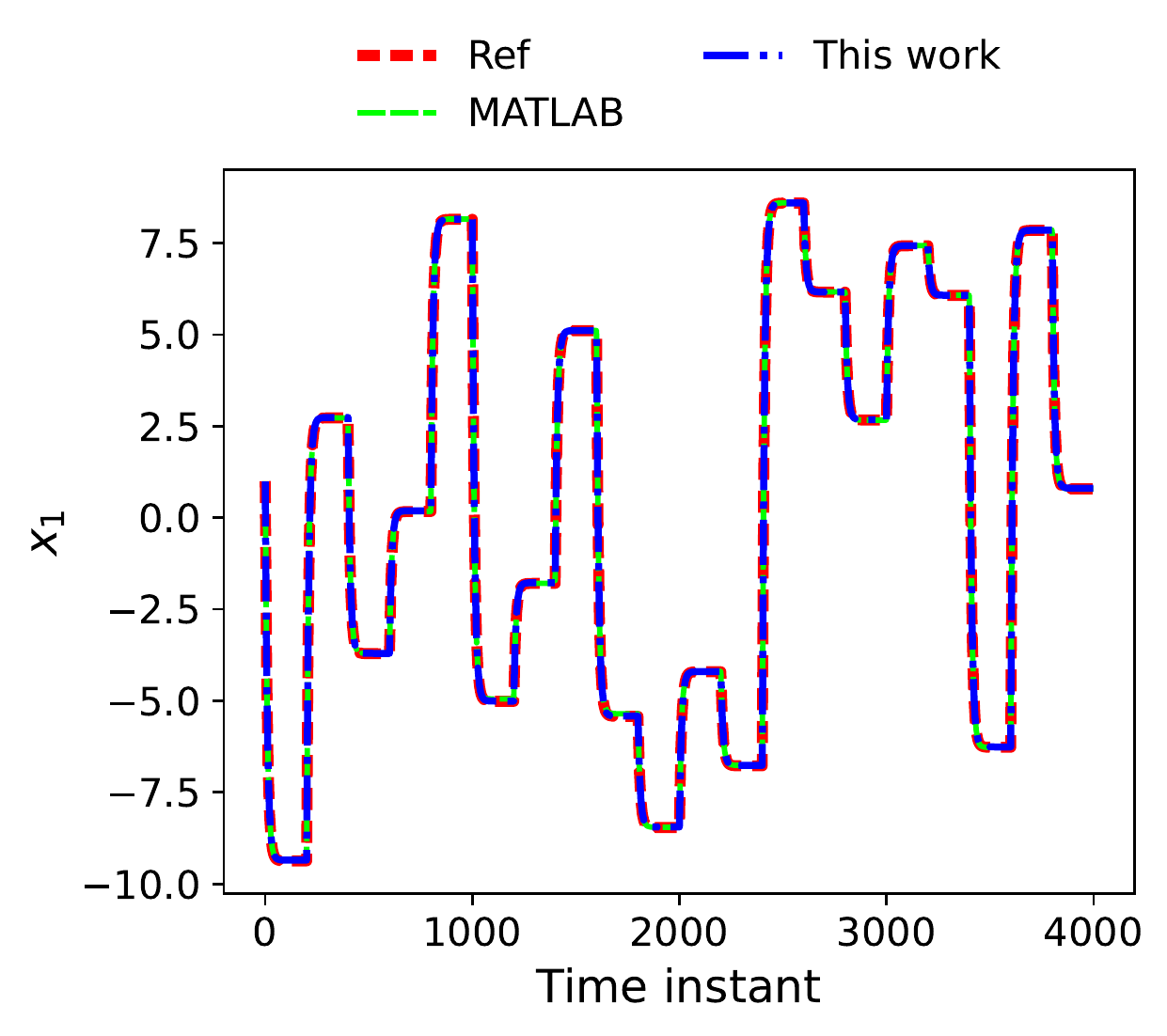}
		\caption{}
	\end{subfigure}
	\begin{subfigure}[]{0.45\textwidth}
		\centering
		\includegraphics[width=\textwidth]{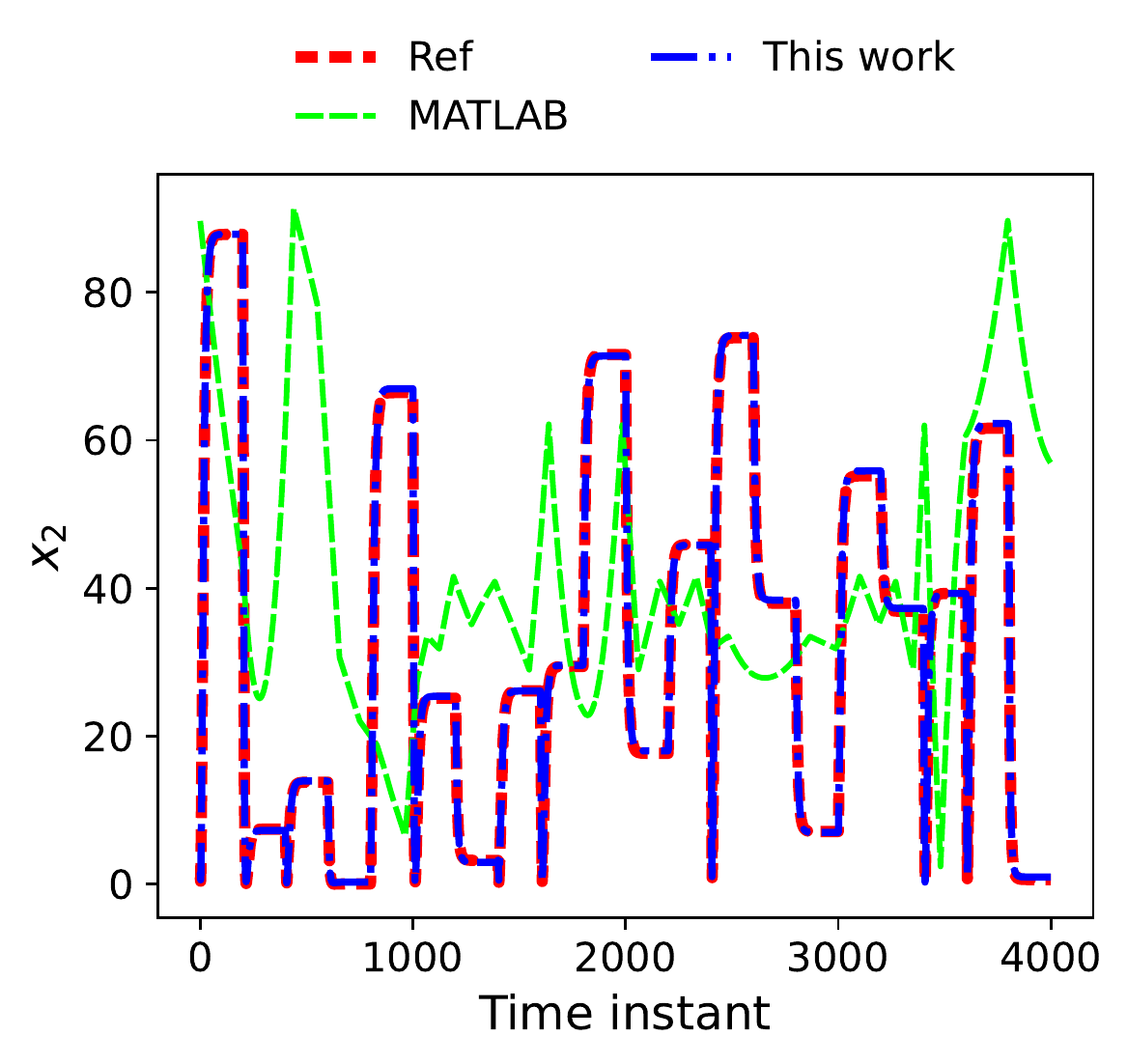}
		\caption{}
	\end{subfigure}
	\begin{subfigure}[]{0.45\textwidth}
		\centering
		\includegraphics[width=\textwidth]{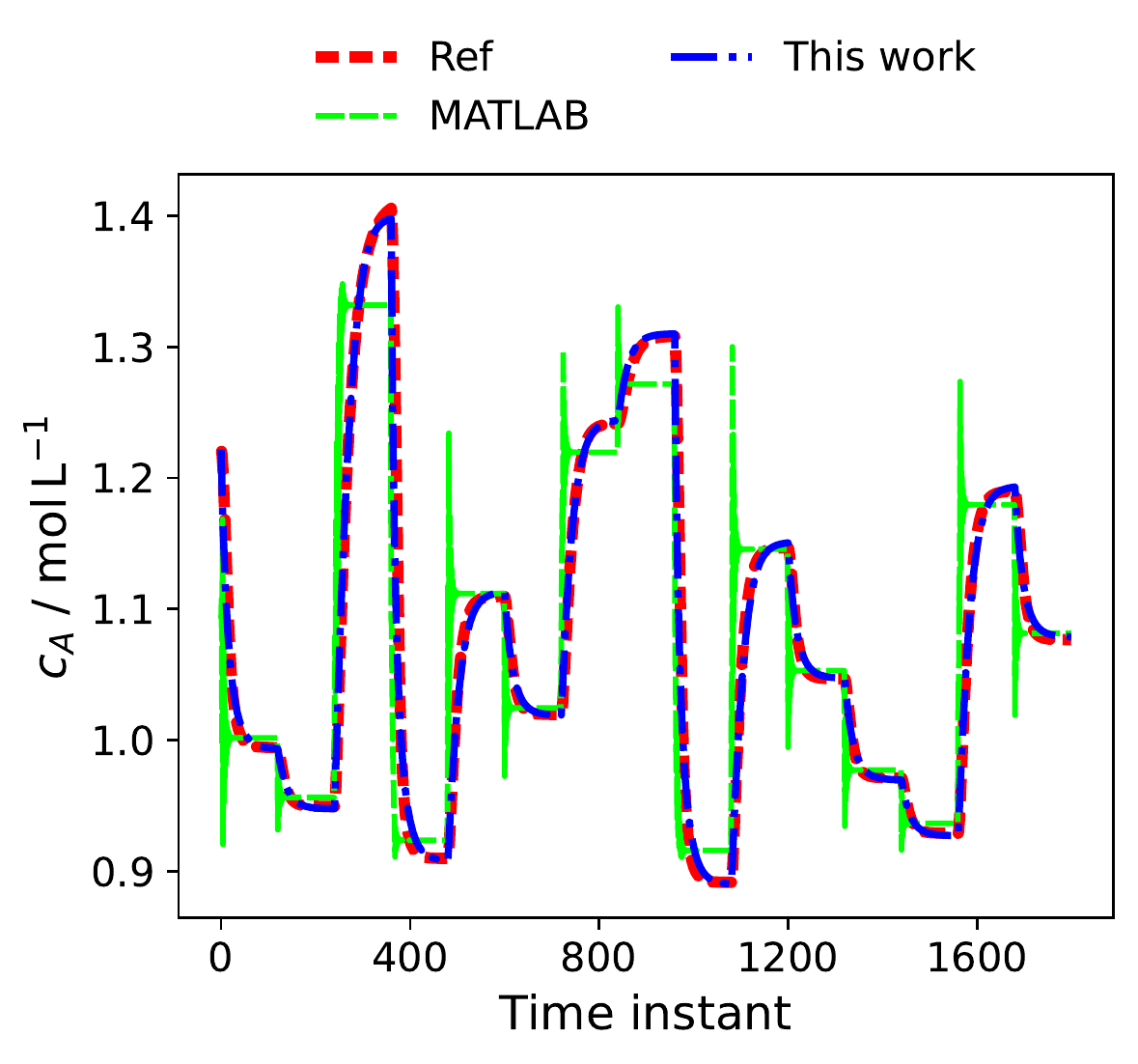}
		\caption{}
	\end{subfigure}
	\begin{subfigure}[]{0.45\textwidth}
		\centering
		\includegraphics[width=\textwidth]{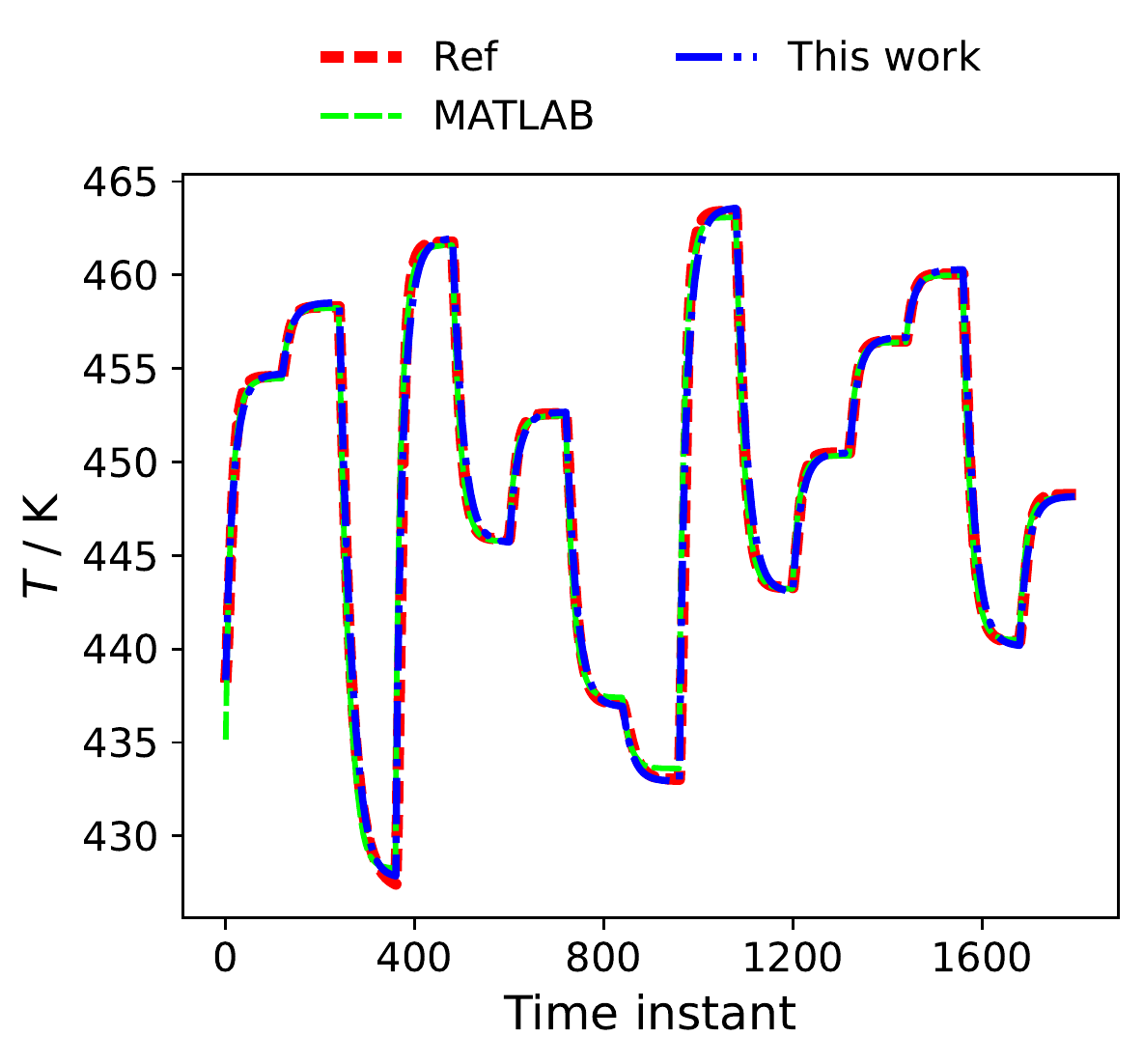}
		\caption{}
	\end{subfigure}
	\begin{subfigure}[]{0.48\textwidth}
		\centering
		\includegraphics[width=\textwidth]{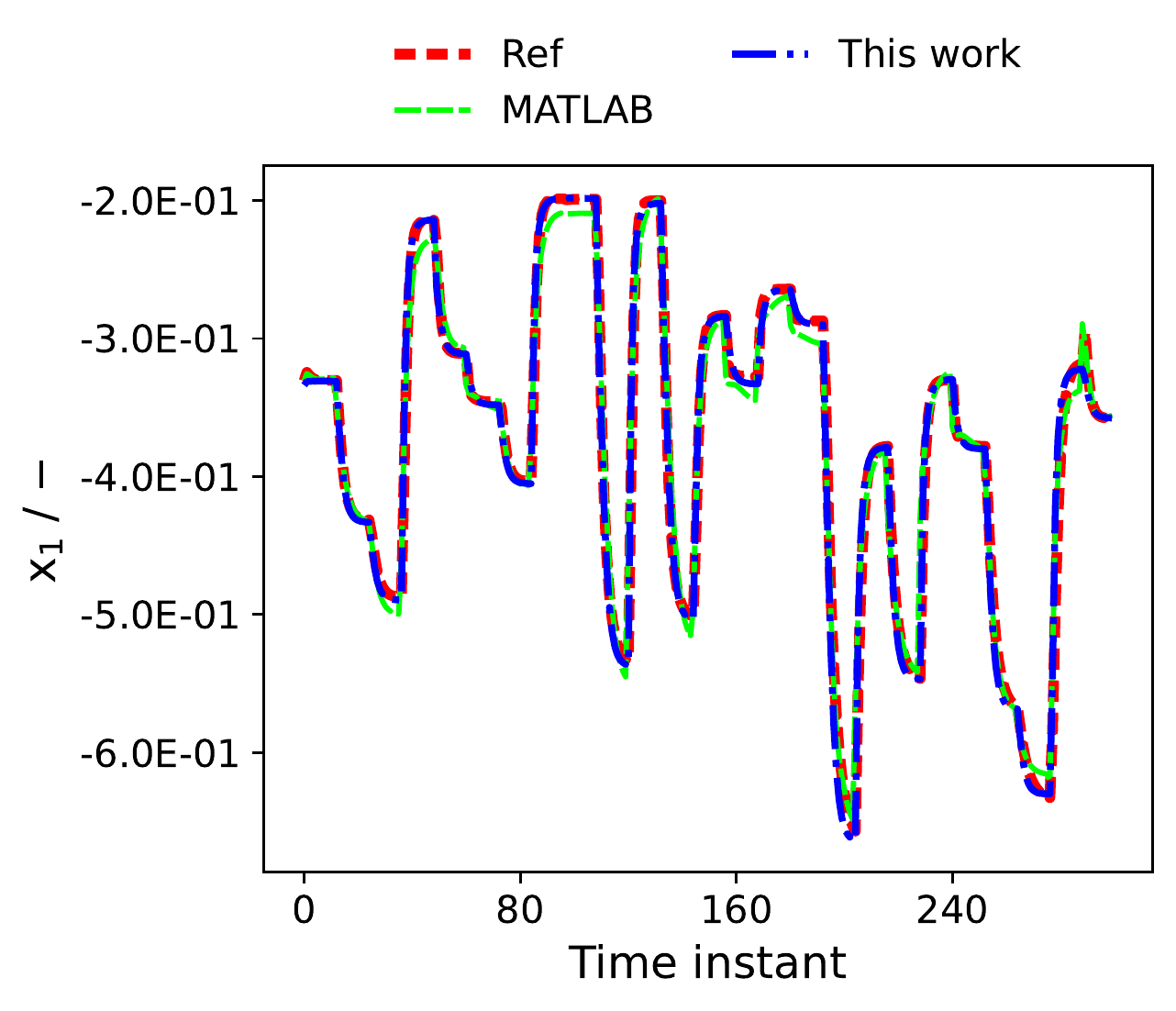}
		\caption{}
	\end{subfigure}
	\begin{subfigure}[]{0.48\textwidth}
		\centering
		\includegraphics[width=\textwidth]{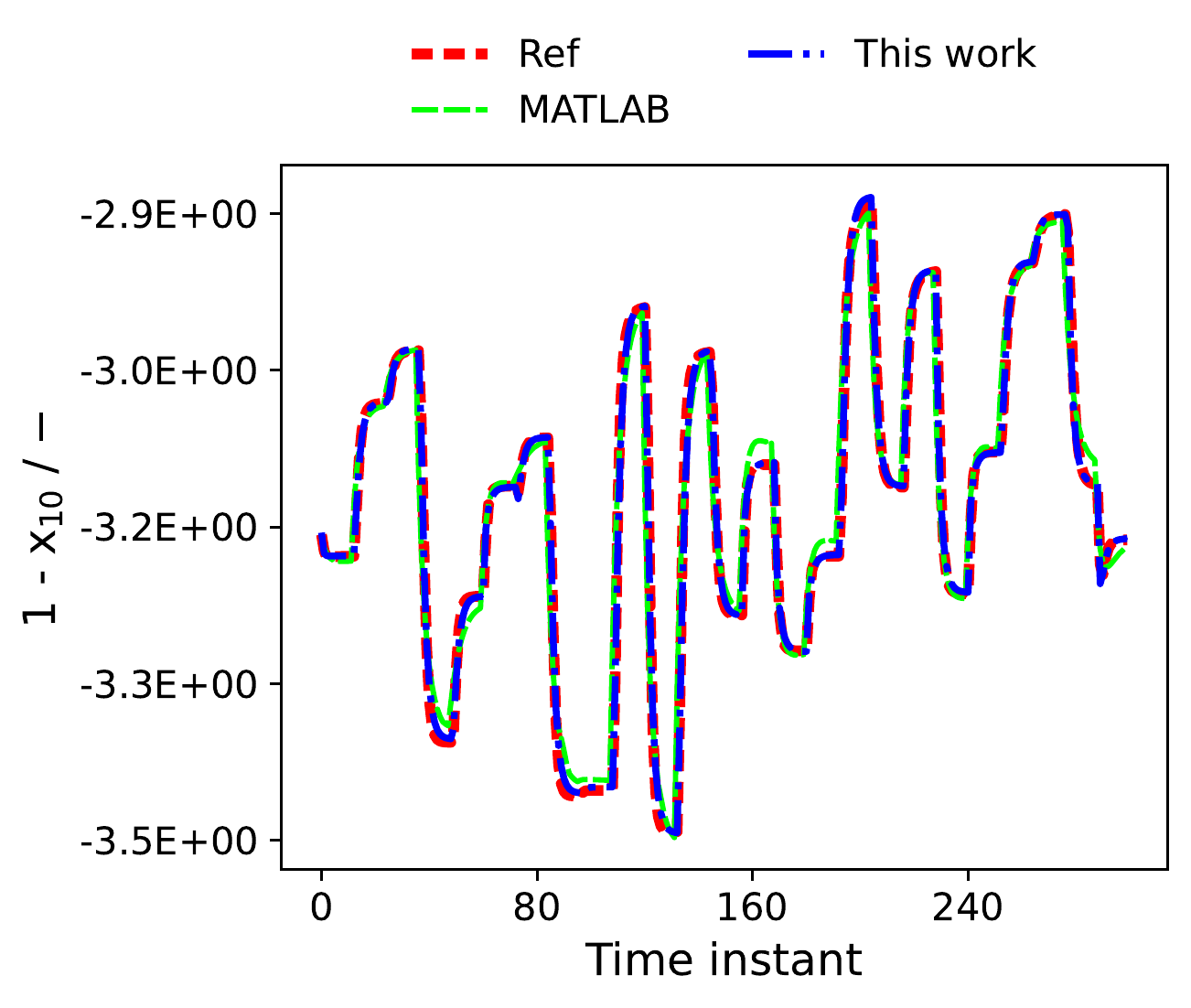}
		\caption{}
	\end{subfigure}
	\caption{Comparison of MIMO Wiener models from MATLAB Toolbox and our deep-learning framework. (a) + (b): Case study 1. (c) + (d): Case study 2. (e) + (f): Case study 3.}
	\label{fig:results_matlab}
\end{figure}

We compare the performance of the MIMO Wiener models identified with our deep-learning framework to MIMO Wiener models generated by the MATLAB 2021a System Identification Toolbox \cite{matlab}.
To this end, we use similar model hyperparameters in both frameworks, \cref{tab:hyperparameters_matlab}.
Further, we use a piecewise-linear nonlinear block with MATLAB, corresponding to ReLU rather than ELU activation since the latter is not provided.
We use the default settings of the MATLAB toolbox.
However, we increase the number of major iterations of the MATLAB algorithm to 100 and use $\lambda = 0.001$ regularization, as we find that these adjustments improve the models considerably.
Both identification frameworks receive the same scaled training data sets.
Notice that the MATLAB toolbox jointly identifies a set of MISO models used as MIMO model, where $n_z$ and $n\ind{n}$ parameterize each MISO model.
Hence, the MATLAB models actually involve $n_y$ times as many parameters and equations as our Wiener-type Koopman models.
For case study 3, we identify a MIMO model only for the target states $x_1$ and $x_{10}$ with the MATLAB toolbox, as the generated full-state MIMO Wiener model performs significantly worse (NMSE = 0.1).
The comparison is visualized in \cref{fig:results_matlab}.
In all case studies, the NMSE of the MIMO Wiener models generated with MATLAB is at least one order of magnitude greater.
Moreover, in case studies 1 and 2, the obtained MATLAB models completely fail to predict one of the states.
The results also underpin that a single MIMO model identified with a Koopman deep-learning approach can provide a more accurate prediction with fewer trainable parameters than a set of MISO models.

\begin{table}[ht]
	\small
	\caption{\small Hyperparameters and prediction error of the MIMO Wiener models: $n_z$ is subspace dimension, $n\ind{n}$ is no.~of hidden neurons of single-hidden-layer encoder/decoder or no.~of linear segments.}
	\label{tab:hyperparameters_matlab}
	\centering
	\begin{tabular}{l l l l l}
		\toprule
		&&& \multicolumn{2}{c}{\textbf{NMSE}} \\ 
		\textbf{MIMO Wiener model} & $\bm{n_z}$ & $\bm{n\ind{n}}$ & \textbf{MATLAB} & \textbf{This work}\\
		\midrule 
		Case Study 1 & 2 & 20  & 0.2 & 3$\cdot 10 ^{-5}$\\[1ex]
		Case Study 2 & 1 & 20  & 0.025 & 0.0019\\[1ex]
		Case Study 3 & 2 & 10  & 0.0024 & 0.0003\\
		\bottomrule
	\end{tabular}
\end{table}

\end{document}